\documentclass{elsarticle}
\usepackage{amsfonts}
\usepackage{amssymb}
\usepackage{amsmath}
\usepackage{moresize}

\makeatletter
\def\ps@pprintTitle{%
 \let\@oddhead\@empty
 \let\@evenhead\@empty
 \def\@oddfoot{}%
 \let\@evenfoot\@oddfoot}
\makeatother

\usepackage{graphicx}

\usepackage{tikz}

\usepackage{xargs}  

\newcommand{\gridLine}[4]{
    \draw[thick, <->, >=latex] (#1, #2) -- (#3, #4);
}

\newcommand{\gridLineSlim}[4]{
    \draw[<->, >=latex] (#1, #2) -- (#3, #4);
}

\def\drawingTreeT{
    \begin{tikzpicture}
    \vertexInCycleTypeOne{0}{0}{0}{0.75}{160}{v1};
    \vertexInCycleTypeOne{0}{0}{0}{0.75}{280}{v11};
    \vertexInCycleTypeOne{0.50}{1.25}{0}{0.75}{180}{v2};
    \vertexInCycleTypeOne{0.50}{1.25}{0}{0.75}{60}{v22};
    \vertexInCycleTypeZero{0.75}{0.75}{0}{120}{v3};
    \vertexInCycleTypeZero{1.5}{0.75}{0}{120}{v4};
    \vertexInCycleTypeOne{2.1}{1.15}{0}{0.75}{105}{v5};
    \vertexInCycleTypeOne{1.75}{0}{0}{0.75}{240}{v6};
    \vertexInCycleTypeOne{3}{1}{0}{0.75}{45}{v7};
    \vertexInCycleTypeOne{3}{1}{0}{0.75}{285}{v77};
    \vertexInCycleTypeZero{2.75}{-0.25}{0}{120}{v8};
    \vertexInCycleTypeOne{3.75}{0}{0}{0.75}{80}{v9};
    \vertexInCycleTypeOne{3.75}{0}{0}{0.75}{320}{v99};
    \vertexInCycleTypeOne{2.75}{-1.25}{0}{0.75}{330}{v0};
    \vertexInCycleTypeOne{2.75}{-1.25}{0}{0.75}{210}{v00};

    \edge{v1}{v3};
    \edge{v2}{v3};
    \edge{v3}{v4};
    \edge{v4}{v5};
    \edge{v4}{v6};
    \edge{v5}{v7};
    \edge{v6}{v8};
    \edge{v8}{v9};
    \edge{v8}{v0};
    \end{tikzpicture}
}

\def\drawingHzero{
    \begin{tikzpicture}
    \vertexInCycleTypeOne{0}{0}{0}{1}{90}{v};
    \vertexInCycleTypeOne{0}{0}{0}{1}{210}{v};
    \vertexInCycleTypeOne{0}{0}{0}{1}{330}{v};
    \end{tikzpicture}
}

\def\drawingExample{
    \begin{tikzpicture}

    \vertex{0}{0}{v3};
    \vertex{0}{1.5}{v1};
    \vertexLabel[above]{v1}{\small $v_1$};
    \vertex{-0.433}{0.75}{v4};
    \vertexLabel[left]{v4}{\small $v_4$};
    \vertex{0.433}{0.75}{v2};
    \vertexLabel[right]{v2}{\small $v_2$};

    \vertex{0.866}{0}{v5};
    \vertexLabel[right]{v5}{\small $v_5$};
    \vertex{-0.866}{0}{v9};
    \vertexLabel[left]{v9}{\small $v_9$};
    \vertex{-0.433}{-0.75}{v8};
    \vertexLabel[below]{v8}{\small $v_8$};
    \vertex{0.433}{-0.75}{v7};
    \vertexLabel[below]{v7}{\small $v_7$};
    \vertex{1.299}{-0.75}{v6};
    \vertexLabel[right]{v6}{\small $v_6$};

    \edge{v3}{v1}
    \edge{v3}{v2}
    \edge{v3}{v4}
    \edge{v1}{v4}
    \edge{v1}{v2}
    \edge{v2}{v4}

    \edge{v3}{v8}
    \edge{v3}{v9}
    \edge{v8}{v9}

    \edge{v3}{v5}
    \edgeLabel[above]{v3}{v5}{\small $v_3$};
    \edge{v3}{v7}
    \edge{v5}{v7}
    \edge{v6}{v5}
    \edge{v6}{v3}
    \edge{v6}{v7}

    \draw[dotted] (3.0,-0.5) -- (5.0,-0.5);
    \draw[dotted] (3.0, 0.0) -- (5.0, 0.0);
    \draw[dotted] (3.0, 0.5) -- (5.0, 0.5);
    \draw[dotted] (3.0, 1.0) -- (5.0, 1.0);
    \draw[dotted] (3.0, 1.5) -- (5.0, 1.5);

    \draw[dotted] (3.0,-0.5) -- (3.0, 1.5);
    \draw[dotted] (3.5,-0.5) -- (3.5, 1.5);
    \draw[dotted] (4.0,-0.5) -- (4.0, 1.5);
    \draw[dotted] (4.5,-0.5) -- (4.5, 1.5);
    \draw[dotted] (5.0,-0.5) -- (5.0, 1.5);

    \gridLine{ 4.0}{ 0.0}{ 4.0}{ 1.0}; 
    \gridLine{ 4.0}{ 0.0}{ 4.0}{-0.5}; 
    \gridLine{ 4.0}{ 1.0}{ 4.0}{ 1.5}; 
    \gridLine{ 3.0}{ 0.0}{ 4.0}{ 0.0}; 
    \gridLine{ 3.0}{ 0.5}{ 4.0}{ 0.5}; 
    \gridLine{ 3.0}{ 1.0}{ 4.0}{ 1.0}; 
    \gridLine{ 5.0}{ 0.0}{ 4.0}{ 0.0}; 
    \gridLine{ 5.0}{ 0.5}{ 4.0}{ 0.5}; 
    \gridLine{ 5.0}{ 1.0}{ 4.0}{ 1.0}; 

    \node [below] at (4.0,-0.5) {\footnotesize $P_{v_6}$};
    \node [above] at (4.0, 1.5) {\footnotesize $P_{v_1}$};
    \node [left]  at (3.0, 0.0) {\footnotesize $P_{v_7}$};
    \node [left]  at (3.0, 0.5) {\footnotesize $P_{v_9}$};
    \node [left]  at (3.0, 1.0) {\footnotesize $P_{v_4}$};
    \node [right] at (5.0, 0.0) {\footnotesize $P_{v_5}$};
    \node [right] at (5.0, 0.5) {\footnotesize $P_{v_8}$};
    \node [right] at (5.0, 1.0) {\footnotesize $P_{v_2}$};
    \node [right] at (3.9, 0.75) {\footnotesize $P_{v_3}$};

    \end{tikzpicture}
}

\def\drawingSplit{
    \begin{tikzpicture}[scale=0.7]

    \draw[dotted] (0,1) -- (4,1);
    \draw[dotted] (0,2) -- (4,2);
    \draw[dotted] (0,3) -- (4,3);
    \draw[dotted] (1,0) -- (1,4);
    \draw[dotted] (4,0) -- (4,4);

    \draw[->, >=latex] (5.25,2) -- (6.25,2);

    \draw[dotted] (8,1) -- (12,1);
    \draw[dotted] (8,2) -- (12,2);
    \draw[dotted] (8,3) -- (12,3);
    \draw[dotted] (9,0) -- (9,4);
    \draw[dotted] (12,0) -- (12,4);

    \node [left]  at (0,1) {\scriptsize $x_{\ell}$};
    \node [left]  at (0,2) {\scriptsize $x_i$};
    \node [left]  at (0,3) {\scriptsize $x_t$};
    \node [below] at (1,0) {\scriptsize $y_j$};
    \node [below] at (4,0) {\scriptsize $y_k$};

    \node [left]  at (8,1) {\scriptsize $x_{\ell}$};
    \node [left]  at (8,2) {\scriptsize $x_i$};
    \node [left]  at (8,3) {\scriptsize $x_t$};
    \node [below] at (9,0) {\scriptsize $y_j$};
    \node [below] at (12,0) {\scriptsize $y_k$};

    \gridLine{1}{1}{1}{3};
    \gridLine{1}{2}{4}{2};
    \gridLine{9}{1}{9}{2};
    \gridLine{9}{2}{9}{3};
    \gridLine{9}{2}{12}{2};

    \node [below] at (1.35,3) {\scriptsize $P_{u}$};
    \node [below] at (9.45,3) {\scriptsize $P_{u_1}$};
    \node [above] at (9.45,1) {\scriptsize $P_{u_2}$};

    \node [above] at (2.5,2) {\scriptsize $P_{v}$};
    \node [above] at (11.5,2) {\scriptsize $P_{v}$};

    \end{tikzpicture}
}

\def\drawingLemmaEndFree{
    \begin{tikzpicture}[scale=0.75]

    \gridLine{0}{2}{1}{2};
    \gridLine{1}{2}{6}{2}; 
    \gridLine{1}{1}{1}{2};
    \gridLine{1}{2}{1}{3};
    \gridLine{2}{1}{2}{2};
    \gridLine{3}{1}{3}{2};
    \gridLine{3}{2}{3}{3};
    \gridLine{4.5}{2}{4.5}{3};

    \node [above] at (6.35,2) {$P_{v}$};
    \node [below] at (2.5,0.9) {$N[v]$};
    \node [below] at (2.5,-0.6) {$(a)$};

    \gridLine{9}{2}{10}{2};
    \gridLine{10}{2}{15}{2}; 
    \gridLine{10}{1}{10}{2};
    \gridLine{10}{2}{10}{3};
    \gridLine{11}{1}{11}{2};
    \gridLine{12}{1}{12}{2};
    \gridLine{12}{2}{12}{3};
    \gridLine{13.5}{2}{13.5}{3};

    \node [above] at (15.35,2) {$P_{v}$};
    \node [below] at (11.5,-0.6) {$(b)$};

    \draw (8.65,1) rectangle (9.35,3.5);
    \draw (9.5,2.5) rectangle (10.5,3.5);
    \draw (11.25,2.5) rectangle (12.75,3.25);
    \draw (13,2.5) rectangle (14.25,3.25);

    \draw (9.65,1.5) rectangle (10.35,0.5);
    \draw (10.6,1.5) rectangle (11.4,0.5);
    \draw (11.65,1.5) rectangle (13.5,0.65);

    \draw (8.45,0.25) rectangle (14.5,3.75);

    \end{tikzpicture}
}

\def\drawingTreeAlgo{
    \begin{tikzpicture}[scale=1.25]
    \vertexInCycleTypeOne{0}{0}{0}{0.75}{160}{v1};
    \vertexInCycleTypeOne{0}{0}{0}{0.75}{280}{v11};
    \vertexInCycleTypeOne{0.50}{1.25}{0}{0.75}{190}{v2};
    \vertexInCycleTypeOne{0.50}{1.25}{0}{0.75}{70}{v22};
    \vertexInCycleTypeZero{0.75}{0.75}{0}{120}{v3};
    \vertexInCycleTypeZero{1.5}{0.75}{0}{120}{v4};
    \vertexInCycleTypeZero{1.125}{0.25}{0}{0}{v10};
    \vertexInCycleTypeOne{2.1}{1.15}{0}{0.75}{105}{v5};
    \vertexInCycleTypeOne{3}{1}{0}{0.75}{45}{v7};
    \vertexInCycleTypeOne{3}{1}{0}{0.75}{285}{v77};
    \vertexInCycleTypeOne{2.25}{0}{0}{0.75}{330}{v0};
    \vertexInCycleTypeOne{2.25}{0}{0}{0.75}{210}{v00};

   \vertexLabel[above]{v1}{\small $1$};
   \vertexLabel[right]{v2}{\small $2$};
   \vertexLabel[left]{v3}{\small $5$};
   \vertexLabel[above]{v4}{\small $6$};
   \node [above] at (2.4,-0.05) {\small $4$};
   \node [above] at (2.9,1) {\small $3$};

    \edge[<-,>=latex]{v1}{v3};
    \edge[<-,>=latex]{v2}{v3};
    \edge[<-,>=latex]{v3}{v10};
    \edge{v4}{v10};
    \edge[<-,>=latex]{v3}{v4};
    \edge{v4}{v5};
    \edge[->,>=latex]{v4}{v0};
    \edge[->,>=latex,thick]{v5}{v7};
    \edge[->,>=latex,thick]{v5}{v0};
    \vertexInCycleTypeZero{1.125}{0.25}{0}{0}{v1010};
    \vertexInCycleTypeZero{1.5}{0.75}{0}{120}{v44};
    \vertexInCycleTypeZero{2.1}{1.15}{0}{0}{v55};
    \node [below] at (2.3,1.1) {\small $R$};

    \end{tikzpicture}
}

\def\drawingTreeCograph{
    \begin{tikzpicture}[scale=0.75]

    \gridLine{7.5}{-1}{7.5}{1};
    \gridLine{7.5}{1}{7.5}{5}; 
    \gridLine{7.5}{5}{7.5}{7};
    \gridLine{6}{1}{7.5}{1};
    \gridLine{9}{1}{7.5}{1};
    \gridLine{9}{2}{7.5}{2};
    \gridLine{6}{3}{7.5}{3};
    \gridLine{9}{3}{7.5}{3};
    \gridLine{6}{4}{7.5}{4};
    \gridLine{9}{4}{7.5}{4};
    \gridLine{6}{5}{7.5}{5};
    \gridLine{9}{5}{7.5}{5};

    \node [left] at (7.5,2.25) {\footnotesize $P_{v_{12}}$};
    \node [above] at (6.75,1) {\footnotesize $P_{v_5}$};
    \node [above] at (6.75,3) {\footnotesize $P_{v_9}$};
    \node [above] at (6.75,4) {\footnotesize $P_{v_7}$};
    \node [above] at (6.75,5) {\footnotesize $P_{v_2}$};

    \node [above] at (8.25,1) {\footnotesize $P_{v_6}$};
    \node [above] at (8.25,2) {\footnotesize $P_{v_{11}}$};
    \node [above] at (8.25,3) {\footnotesize $P_{v_{10}}$};
    \node [above] at (8.25,4) {\footnotesize $P_{v_8}$};
    \node [above] at (8.25,5) {\footnotesize $P_{v_3}$};

    \node [right] at (7.5,0) {\footnotesize $P_{v_4}$};
    \node [right] at (7.5,6) {\footnotesize $P_{v_1}$};

    \node [above] at (7.5,7.5) {\footnotesize Representation};
    \node [above] at (1,7.5) {\footnotesize $G$};

    \vertex{1}{-0.5}{v11};
    \vertex{0}{1.25}{v9};
    \vertex{2}{1.25}{v10};
    \vertex{0}{2.25}{v7};
    \vertex{2}{2.25}{v8};
    \vertex{0}{3.25}{v5};
    \vertex{2}{3.25}{v6};
    \vertex{1}{4.35}{v4};
    \vertex{0}{5.45}{v2};
    \vertex{2}{5.45}{v3};
    \vertex{1}{6.55}{v1};
    \vertex{3.5}{2.75}{v12};

    \edge{v1}{v12};
    \edge{v2}{v12};
    \edge{v3}{v12};
    \edge{v4}{v12};
    \edge{v5}{v12};
    \edge{v6}{v12};
    \edge{v7}{v12};
    \edge{v8}{v12};
    \edge{v9}{v12};
    \edge{v10}{v12};
    \edge{v11}{v12};

    \edge{v9}{v10};
    \edge{v7}{v8};
    \edge{v5}{v6};
    \edge{v4}{v5};
    \edge{v4}{v6};
    \edge{v1}{v2};
    \edge{v1}{v3};
    \edge{v2}{v3};

    \vertexLabel[right]{v12}{\footnotesize $v_{12}$};
    \vertexLabel[left]{v1}{\footnotesize $v_1$};
    \vertexLabel[left]{v2}{\footnotesize $v_2$};
    \vertexLabel[left]{v4}{\footnotesize $v_4$};
    \vertexLabel[left]{v5}{\footnotesize $v_5$};
    \vertexLabel[left]{v7}{\footnotesize $v_7$};
    \vertexLabel[left]{v9}{\footnotesize $v_9$};
    \vertexLabel[left]{v11}{\footnotesize $v_{11}$};
    \vertexLabel[right]{v3}{\footnotesize $v_{3}$};
    \vertexLabel[right]{v6}{\footnotesize $v_{6}$};
    \vertexLabel[right]{v8}{\footnotesize $v_{8}$};
    \vertexLabel[right]{v10}{\footnotesize $v_{10}$};

    \end{tikzpicture}
}

\def\drawingSpider{
    \begin{tikzpicture}[scale=0.75]

    \gridLine{6.5}{2}{8}{2};
    \gridLine{9.5}{2}{8}{2};
    \gridLine{9.5}{2}{11}{2};
    \gridLine{12.5}{2}{11}{2};
    \gridLine{9.5}{-1}{9.5}{0.5};
    \gridLine{9.5}{2}{9.5}{0.5};
    \gridLine{9.5}{2}{9.5}{3.5};
    \gridLine{9.5}{5}{9.5}{3.5};

    \node [left] at (9.5,1.25) {\footnotesize $P_{c_3}$};
    \node [left] at (9.5,-0.25) {\footnotesize $P_{s_3}$};

    \node [above] at (8.75,2) {\footnotesize $P_{c_4}$};
    \node [above] at (7.25,2) {\footnotesize $P_{s_4}$};

    \node [right] at (9.5,2.75) {\footnotesize $P_{c_1}$};
    \node [right] at (9.5,4.25) {\footnotesize $P_{s_1}$};

    \node [below] at (10.25,2) {\footnotesize $P_{c_2}$};
    \node [below] at (11.75,2) {\footnotesize $P_{s_2}$};

    \node [above] at (9.5,5.5) {\footnotesize Representation};
    \node [above] at (2,5.5) {\footnotesize $|C| = 4$ and $R = \emptyset$};

    \vertex{1}{1}{c4};
    \vertex{3}{1}{c3};
    \vertex{1}{3}{c1};
    \vertex{3}{3}{c2};
    \vertex{0}{0}{s4};
    \vertex{0}{4}{s1};
    \vertex{4}{0}{s3};
    \vertex{4}{4}{s2};

    \edge{c1}{c2};
    \edge{c1}{c3};
    \edge{c1}{c4};
    \edge{c1}{s1};
    \edge{c2}{c3};
    \edge{c2}{c4};
    \edge{c2}{s2};
    \edge{c3}{c4};
    \edge{c3}{s3};
    \edge{c4}{s4};

    \vertexLabel[above]{c1}{\footnotesize $c_1$};
    \vertexLabel[above]{c2}{\footnotesize $c_2$};
    \vertexLabel[below]{c3}{\footnotesize $c_3$};
    \vertexLabel[below]{c4}{\footnotesize $c_4$};
    \vertexLabel[left]{s1}{\footnotesize $s_1$};
    \vertexLabel[left]{s4}{\footnotesize $s_4$};
    \vertexLabel[right]{s2}{\footnotesize $s_2$};
    \vertexLabel[right]{s3}{\footnotesize $s_3$};

    \end{tikzpicture}
}

\def\drawingFatSpidersRep{
    \begin{tikzpicture}[scale=0.375]

    \gridLineSlim{0}{-8}{2}{-8};
    \gridLineSlim{2}{-8}{4}{-8};
    \gridLineSlim{4}{-8}{6}{-8};
    \gridLineSlim{6}{-8}{8}{-8};
    \gridLineSlim{4}{-8}{4}{-6};
    \gridLineSlim{4}{-6}{4}{-4};
    \gridLineSlim{4}{-6}{6}{-6};
    \gridLineSlim{4}{-8}{4}{-10};
    \gridLineSlim{4}{-10}{6}{-10};

    \gridLineSlim{10}{-8}{12}{-8};
    \gridLineSlim{12}{-8}{14}{-8};
    \gridLineSlim{14}{-8}{16}{-8};
    \gridLineSlim{16}{-8}{18}{-8};
    \gridLineSlim{14}{-8}{14}{-6};
    \gridLineSlim{14}{-6}{14}{-4};
    \gridLineSlim{14}{-7}{16}{-7};
    \gridLineSlim{14}{-8}{14}{-10};
    \gridLineSlim{14}{-10}{16}{-10};

    \gridLineSlim{22}{-8}{24}{-8};
    \gridLineSlim{22}{-6}{24}{-6};
    \gridLineSlim{22}{-8}{22}{-6};
    \gridLineSlim{24}{-8}{24}{-6};
    \gridLineSlim{24}{-7}{26}{-7};

        \vertex[white]{28.5}{-8}{s2};

    \end{tikzpicture}
}

\newcommand{\vertex}[4][black]{
    \draw[#1, fill=#1, inner sep=0pt] (#2, #3) circle (0.075) node(#4){};
}

\newcommand{\vertexLabel}[3][above]{
    \path (#2) node[#1]{#3};
}

\newcommand{\edge}[3][]{
	\draw[#1] (#2) -- (#3);
}

\newcommand{\edgeLabel}[4][above]{
	\path (#2) -- (#3) node[midway, #1]{#4};
}

\newcommand{\vertexRegularPolygon}[7][]{
	\vertex[#1]
	{ { (#2) + (#4) * -cos(deg(6.283*(#6 - 1)/(#5)-(1.571+3.142/(#5)))) } }
	{ { (#3) + (#4) * sin(deg(6.283*(#6 - 1)/(#5)-(1.571+3.142/(#5)))) } }
	{#7};
}

\newcommand{\fatSpider}[6][]{
	\foreach \fatSpiderI in {1,...,#5} {
		\vertexRegularPolygon[#1]
						{#2}{#3}
						{#4}
						{#5}
						{\fatSpiderI}
						{#6c\fatSpiderI}
						{};

		\vertexRegularPolygon[#1]
						{#2}{#3}
						{2*(#4)}
						{#5}
						{\fatSpiderI}
						{#6s\fatSpiderI}
						{};
	}

	\foreach \fatSpiderI in {2,...,#5} {
		\pgfmathtruncatemacro{\fatSpiderIminusOne}{\fatSpiderI - 1};

		\foreach \fatSpiderJ in {1,...,\fatSpiderIminusOne} {
			\edge[#1]{#6c\fatSpiderI}{#6c\fatSpiderJ};
		}
	}

	\foreach \fatSpiderI in {1,...,#5} {
		\foreach \fatSpiderJ in {1,...,#5} {
			\ifthenelse
			{\equal{\fatSpiderI}{\fatSpiderJ}}
			{}
			{\edge[#1]{#6s\fatSpiderI}{#6c\fatSpiderJ};};
		}
	}
}

\newcommand{\vertexInCycleTypeZero}[6][black]{
    \vertex[#1]
    { { (#2) + (#4) * cos(#5) } } 
    { { (#3) + (#4) * sin(#5) } } 
    {#6}; 
}

\newcommand{\vertexInCycleTypeOne}[7][black]{
    \vertexInCycleTypeZero[#1]{#2}{#3}{#4}{#6}{#7};

    \vertex[#1]
    { { (#2)+(#4)*cos(#6) + (sqrt(3.0)/3.0*(#5)) * cos((#6)+30.0) } } 
    { { (#3)+(#4)*sin(#6) + (sqrt(3.0)/3.0*(#5)) * sin((#6)+30.0) } } 
    {#7a}; 
    \vertex[#1]
    { { (#2) + ((#4)+(#5)) * cos(#6) } } 
    { { (#3) + ((#4)+(#5)) * sin(#6) } } 
    {#7b}; 
    \vertex[#1]
    { { (#2)+(#4)*cos(#6) + (sqrt(3.0)/3.0*(#5)) * cos((#6)-30.0) } } 
    { { (#3)+(#4)*sin(#6) + (sqrt(3.0)/3.0*(#5)) * sin((#6)-30.0) } } 
    {#7c}; 

    \edge[color=#1]{#7}{#7a};
    \edge[color=#1]{#7}{#7b};
    \edge[color=#1]{#7}{#7c};
    \edge[color=#1]{#7a}{#7b};
    \edge[color=#1]{#7a}{#7c};
    \edge[color=#1]{#7b}{#7c};
}

\newcommand{\thin}{\mbox{thin}}
\newcommand{\thick}{\mbox{thick}}
\newcommand{\conB}{contact $B_0$-VPG}
\newtheorem{thm}{Theorem}
\newtheorem{lem}[thm]{Lemma}
\newtheorem{rem}[thm]{Remark}
\newproof{pf}{Proof}
\newtheorem{cor}[thm]{Corollary}
\newtheorem{fact}[thm]{Fact}

\begin{document}


    \begin{frontmatter}
        \title{On some special classes of contact $B_0$-VPG graphs}
        \author[buenos,ICC]{Flavia Bonomo-Braberman}
        \ead{fbonomo@dc.uba.ar}
        \author[laplata,CEMaLP]{Mar\'{\i}a P\'{\i}a Mazzoleni}
        \ead{pia@mate.unlp.edu.ar}
            \author[buenos,ICC]{Mariano Leonardo Rean}
            \ead{marianorean@gmail.com}
            \author[suiza]{Bernard Ries}
            \ead{bernard.ries@unifr.ch}
            \address[buenos]{Universidad de Buenos Aires. Facultad de Ciencias
Exactas y Naturales. Departamento de Computaci\'on. Buenos Aires,
Argentina.}
        \address[laplata]{Departamento de Matem\'atica,
            Universidad Nacional de La Plata, La Plata,
            Argentina.}
        \address[suiza]{University of Fribourg,
            Fribourg, Switzerland.}
    \address[ICC]{CONICET-Universidad de Buenos Aires. Instituto de
Investigaci\'on en Ciencias de la Computaci\'on (ICC). Buenos
Aires, Argentina.}
            \address[CEMaLP]{Centro de Matem\'atica Aplicada de La Plata (CEMaLP), La Plata, Argentina.}
            %


\begin{abstract}
A graph $G$ is a $B_0$-VPG graph if one can associate a horizontal
or vertical path on a rectangular grid with each vertex such that
two vertices are adjacent if and only if the corresponding paths
intersect in at least one  grid-point. A graph $G$ is a
\emph{contact $B_0$-VPG graph} if it is a $B_0$-VPG graph
admitting a representation with no one-point paths, no two paths
crossing, and no two paths sharing an edge of the grid. In this
paper, we present a minimal forbidden induced subgraph
characterisation of contact $B_0$-VPG graphs within four special
graph classes: chordal graphs, tree-cographs, $P_4$-tidy graphs
and $P_5$-free graphs. Moreover, we present a polynomial-time
algorithm for recognising chordal contact $B_0$-VPG graphs.

{\em Keywords: contact $B_0$-VPG graph, chordal graph,
tree-cograph, $P_4$-tidy graph, $P_5$-free graph.}

\end{abstract}
\end{frontmatter}

\section{Introduction}

Golumbic et al. introduced in~\cite{asinowski} the concept of
\textit{vertex intersection graphs of paths in a grid} (referred
to as \textit{VPG graphs}). An undirected graph $G=(V,E)$ is
called a VPG graph if one can associate a path in a rectangular
grid with each vertex such that two vertices are adjacent if and
only if the corresponding paths intersect in at least one
grid-point. In the seminal paper on VPG graphs it was shown that
this class is equivalent to the earlier defined class of string
graphs~\cite{E-E-T-intersec}.

Under the perspective of paths in grids, a particular attention
was paid to the case where the paths have a limited number of
bends. An undirected graph $G=(V,E)$ is then called a
\textit{$B_k$-VPG graph}, for some integer $k\geq 0$, if one can
associate a path with at most $k$ bends in a rectangular grid with
each vertex such that two vertices are adjacent if and only if the
corresponding paths intersect in at least one grid-point.
Recognition of VPG graphs is NP-complete by the equivalence with
string graphs. Moreover $B_k$-VPG recognition is NP-complete for
all $k$~\cite{chaplick12}.

Since their introduction, $B_k$-VPG graphs have been studied by
many researchers and the community of people working on these
graph classes or related ones is still growing (see for
instance~\cite{ABM-VPG-gc, asinowski, chaplick, chaplick13,
cohen1, cohen2, Felsner, Golumbic-Ries}).


In this paper, we are interested in a subclass of $B_k$-VPG graphs
called \textit{contact $B_k$-VPG}. A \textit{contact $B_k$-VPG
representation} of $G$ is a VPG representation in which each path
has length at least one, at most $k$ bends, and intersecting paths
neither cross each other nor share an edge of the grid. A graph is
a \textit{contact $B_k$-VPG graph} if it has a contact $B_k$-VPG
representation. Here, we will focus on the special case when
$k=0$, i.e. each path is a horizontal or vertical path in the
grid.

Contact graphs in general (graphs where vertices represent
geometric objects which are allowed to touch but not to cross each
other, a natural model arising from real physical objects) have
been considered in the past (see for
instance~\cite{castro,Fraysseix1,Hlineny,Hlineny1}). In
particular, for intersection models of lines in the plane, it is
often the case that three lines intersecting at a same point is
not allowed, but we do not impose such a restriction.

As for many graph classes having not many known full
characterisations (for example, a complete list of minimal
forbidden induced subgraphs is not known), their characterisation
within well studied graphs classes or with respect to graph
parameters was investigated. In the case of contact $B_k$-VPG
graphs, it was shown in~\cite{Fraysseix2} that every planar
bipartite graph is a contact $B_0$-VPG graph. Later,
in~\cite{chaplick13}, the authors show that every triangle-free
planar graph is a contact $B_1$-VPG graph. In a recent paper
(see~\cite{D-G-M-R-cpg}), contact $B_k$-VPG graphs have been
investigated from a structural point of view and it was for
instance shown that they do not contain cliques of size $7$ and
they always contain a vertex of degree at most $6$. Moreover, it
was shown that they are $6$-colourable. Regarding contact
$B_0$-VPG graphs, it was shown that they are $4$-colourable.
Furthermore, $3$-colouring and the recognition problem were shown
to be NP-complete.

In this paper, our goal is to get a better understanding and
knowledge of the underlying structure of contact $B_0$-VPG graphs.
Even though classical graph problems may be difficult to solve in
this graph class (see for instance~\cite{D-G-M-R-cpg}), a better
knowledge of its structural properties may lead to good
approximation algorithms for these problems. We will consider the
following four special graph classes: chordal graphs,
tree-cographs, $P_4$-tidy graphs and $P_5$-free graphs, and we
will characterise those graphs from these families that are
contact $B_0$-VPG. Moreover, we will present a polynomial-time
algorithm for recognising chordal contact $B_0$-VPG graphs based
on our characterisation. For the other graph classes considered
here, the characterisation immediately yields a polynomial-time
recognition algorithm.

A preliminary version of this paper appears in~\cite{BMRR-ISCO18}.

\section{Preliminaries}\label{sec:definitions}

For concepts and notations not defined here we refer the reader
to~\cite{Bondy-Murty}. All graphs in this paper are simple (i.e.,
without loops or multiple edges). Let $G=(V,E)$ be a graph. If
$u,v\in V$ and $uv\notin E$, $uv$ is called a \emph{nonedge} of
$G$. We write $G-v$ for the subgraph obtained by deleting a vertex
$v$ and all the edges incident to $v$. Similarly, we write $G-e$
for the subgraph obtained by deleting an edge $e$ without deleting
its endpoints.

For each vertex $v$ of $G$, $N_G(v)$ denotes the
\emph{neighbourhood} of $v$ in $G$ and $N_G[v]$ denotes the
\emph{closed neighbourhood}, i.e. $N_G(v)\cup \{v\}$. For a set
$A\subseteq V$, we denote by $N(A)$ the set of vertices having a
neighbour in $A$, and by $N[A]$ the set of vertices belonging to
$A$ or having a neighbour in $A$. Two vertices $v$ and $w$ of $G$
are \emph{false twins} (resp. \emph{true twins}) if
$N_G(v)=N_G(w)$ (resp.\ $N_G[v]=N_G[w]$).

Given a subset $A\subseteq V$, $G[A]$ stands for  \emph{the
subgraph of $G$ induced by $A$}, and  $G \setminus A$ denotes the
\emph{induced subgraph} $G[V\setminus A]$.  We say that a vertex
$v\in V \setminus A$ is \emph{complete to} $A$ if $v$ is adjacent
to every vertex of $A$, and that $v$ is \emph{anticomplete to} $A$
if $v$ has no neighbour in $A$. Similarly, we say that two
disjoint sets $A,B\subset V$ are \emph{complete} (resp.
\emph{anticomplete}) to each other if every vertex in $A$ is
complete (resp. anticomplete) to $B$.

A \emph{clique} is a set of pairwise adjacent vertices. A vertex $v$ is \emph{simplicial}, if $N_G(v$) is a clique. A \emph{stable set} is a set of vertices no two of which are adjacent. A \emph{complete graph} is a graph such that all its vertices are adjacent to each other, i.e. a graph induced by a clique. The \emph{complete graph} on $n$ vertices is denoted by $K_n$. In particular, $K_3$ is called a \emph{triangle}. \emph{$K_4^-$} stands for the graph obtained from $K_4$ by deleting exactly one edge.

The \emph{complement graph} of $G=(V,E)$ is the graph
$\overline{G}=(V,\overline{E})$ such that $\overline{E}=\{uv|\
uv\not\in E\}$. Let $G_1=(V_1,E_1)$ and $G_2=(V_2,E_2)$ be two
graphs. The \emph{disjoint union of $G_1$ and $G_2$}, denoted by
$G_1\cup G_2$, is the graph whose vertex set is $V_1\cup V_2$ and
whose edge set is $E_1\cup E_2$. The \emph{join of $G_1$ and
$G_2$}, denoted by $G_1\vee G_2$, is the graph obtained by first
taking the disjoint union of $G_1$ and $G_2$ and then making $V_1$
and $V_2$ complete to each other. Notice that $\overline{G_1 \cup
G_2} = \overline{G_1} \vee \overline{G_2}$.

Given a graph $H$, we say that $G$ \emph{contains no induced $H$}, if $G$ contains no induced subgraph isomorphic to $H$. If $\mathcal H$ is a family of graphs, $G$ is said to be \emph{$\mathcal H$-free} if $G$ contains no induced subgraph isomorphic to some graph belonging to $\mathcal H$.

Let $\mathcal G$ be a class of graphs. A graph belonging to $\mathcal G$ is called a \emph{$\mathcal G$-graph}. If $G\in\mathcal G$ implies that every induced subgraph of $G$ is a $\mathcal G$-graph, $\mathcal G$ is said to be \emph{hereditary}. If $\mathcal G$ is a hereditary class, a graph $H$ is a \emph{minimal forbidden induced subgraph of $\mathcal G$}, or more briefly, \emph{minimally non-$\mathcal G$}, if $H$ does not belong to $\mathcal G$ but every proper induced subgraph of $H$ is a $\mathcal G$-graph.

A \emph{path} is a sequence of vertices $v_1,\ldots,v_k$ such that
$v_i$ is adjacent to $v_{i+1}$, for $i=1,\ldots,k-1$. The vertices
$v_2,\ldots,v_{k-1}$ are called \emph{internal vertices} of the
path. If there is no edge $v_iv_j$ such that $\vert i-j\vert\geq
2$, the path is said to be \emph{chordless} or \emph{induced}. A
\emph{cycle} $C$ is a sequence of vertices $v_1,\ldots,v_k$ such
that $v_i$ is adjacent to $v_{i+1}$ for $i=1,\ldots,k$, where
indices are taken modulo $k$. If there is no edge $v_iv_j$ such
that $\vert i-j\vert\geq 2$, $C$ is said to be \emph{chordless} or
\emph{induced}. The induced path (resp. induced cycle) on $n$
vertices is denoted $P_n$ (resp. $C_n$). A graph is called
\emph{chordal} if it does not contain any chordless cycle of
length at least four. A \emph{block graph} is a chordal graph
which is $K_4^-$-free.

A graph is \emph{bipartite}, if its vertex set can be partitioned into two stable sets. If, in addition, the two stable sets are complete to each other, the graph is called \emph{complete bipartite}. $K_{n,m}$ stands for the complete bipartite graph whose vertex set can be partitioned into two stable sets $V_1,V_2$ such that $|V_1|=n$ and $|V_2|=m$.

A graph $G$ is \emph{connected}, if for each pair of vertices $u,v$ there exists a path from $u$ to $v$. A \emph{tree} is a connected graph with no induced cycle. Given a connected graph $G=(V,E)$, the \emph{distance between two vertices $u,v\in V$}, denoted by $d_G(u,v)$, is the number of edges of a shortest path from $u$ to $v$. The \emph{diameter} of $G$ is the maximum distance between two vertices.

\begin{figure}
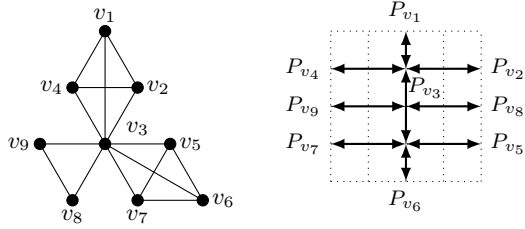

    \begin{center}
    \drawingExample
        \caption{A graph $G$ and a contact $B_0$-VPG
            representation of it.} \label{fig:ejemplo}
    \end{center}
\end{figure}

An undirected graph $G=(V,E)$ is called a \textit{$B_k$-VPG
graph}, for some integer $k\geq 0$, if one can associate a path
with at most $k$ \textit{bends} (a bend is a $90$ degrees turn of
a path at a grid-point) on a rectangular grid with each vertex
such that two vertices are adjacent if and only if the
corresponding paths intersect in at least one  grid-point. Such a
representation is called a \emph{$B_k$-VPG representation}. The
horizontal grid lines will be referred to as \textit{rows} and
denoted by $x_0,x_1,\ldots$ and the vertical grid lines will be
referred to as \textit{columns} and denoted by $y_0,y_1,\ldots$.
We are interested in a subclass of $B_0$-VPG graphs called contact
$B_0$-VPG. A \textit{contact $B_0$-VPG representation}
$\mathcal{R}(G)$ of $G$ is a $B_0$-VPG representation in which
each path in the representation is either a horizontal path or a
vertical path on the grid, with length at least one (the length is
the number of grid-points minus one), such that two vertices are
adjacent if and only if the corresponding paths intersect in at
least one grid-point without crossing each other and without
sharing an edge of the grid. A graph is a \textit{contact
$B_0$-VPG graph} if it has a contact $B_0$-VPG representation. For
every vertex $v$, we denote by $P_v$ the corresponding path in
$\mathcal{R}(G)$ (see Figure~\ref{fig:ejemplo}). Consider a clique
$K$ in $G$. A path $P_v$ representing a vertex $v\in K$ is called
a \emph{path of the clique $K$}.

\

Let us start with an easy but very helpful lemma.

\begin{lem}
    \label{o-clique}
    Let $G$ be a contact $B_0$-VPG graph. Then the size of a biggest clique in $G$ is at most 4, i.e. $G$ is $K_5$-free.
\end{lem}

\begin{pf}
  Given two adjacent vertices in $G$, the intersection of their paths in any contact $B_0$-VPG representation is exactly one grid point. Moreover, it is easy to see that all paths corresponding to vertices in a clique of $G$ must intersect in the same grid point. Assume there is a clique $K$ of size 5 in $G$ and let $P$ be the point of intersection of the corresponding paths in the grid. At least two of the paths must be in the same row or the same column, and contain at least one grid edge intersecting $P$ (a path cannot be only a grid point), a
  contradiction.\qed
\end{pf}

\begin{rem}
    \label{o-Ck} Let $G$ be a $K_4^-$-free graph containing an
    induced cycle $C$ of at least 4 vertices. Then no vertex is
    adjacent to 3 consecutive vertices of $C$.
\end{rem}

Let $G$ be a contact $B_0$-VPG graph, and $K$ be a clique in $G$. A
vertex $v$ is called an \textit{end} in a contact $B_0$-VPG
representation of $K$ if the grid point representing the
intersection of the paths of $K$ corresponds to an endpoint of $P_v$.

\begin{rem}
    \label{o-k4} Let $G$ be a contact $B_0$-VPG graph, and $K$ be a
    clique in $G$ of size four. Then, every vertex in $K$ is an end in
    any contact $B_0$-VPG representation of $K$.
\end{rem}

\begin{lem}
    \label{lem:c4} In any contact $B_0$-VPG representation
    of $C_4$,  the union of the paths representing vertices in $C$ must
    enclose a rectangle of the grid.
\end{lem}

\begin{pf}
  Consider a $B_0$-VPG representation of $C_4$.
  At least two vertices, say $a$ and $b$, in $C$ have the same direction.
  We can assume that $P_a$ and $P_b$ are both vertical.
  If $a$ and $b$ are adjacent, then the corresponding paths intersect in a row $x_i$ of the grid.
  One of them, say $P_a$, is above $x_i$ and the other is below $x_i$.
  Let $c$ be the vertex adjacent to $a$ and non adjacent to $b$.
  Clearly, the path $P_c$ representing $c$ must be also above $x_i$.
  Similarly, the path representing the vertex $d$ adjacent to $b$ and non adjacent to $a$ must be below $x_i$.
  But then it is impossible for $P_c$ and $P_d$ to intersect.
  Therefore, $a$ and $b$ are non adjacent. Now, it is clear that $P_c$ and $P_d$ must be both horizontal,
  otherwise we could repeat the previous argument. If $P_a$ and $P_b$ lie in columns $y_i$ and $y_j$, then
  $P_c$ and $P_d$ must contain all points of the grid between $y_i$ and $y_j$ in their respective columns, say $x_k$ and $x_l$. Then, these paths enclose the rectangle limited by rows $y_i$, $y_j$ and columns $x_k$, $x_l$.
\qed   \end{pf}

\begin{figure}
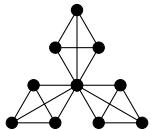

    \begin{center}
\drawingHzero
        \caption{The graph $H_0$.} \label{fig:H0}
    \end{center}
\end{figure}

\

In what follows, we give a set of graphs that are not contact
$B_0$-VPG graphs. We will use this result later to obtain our
characterisations.  Let $H_0$ denote the graph composed of three
$K_4$'s that share a common vertex and such that there are no
further edges (see Figure~\ref{fig:H0}).

\begin{lem}
    \label{lem:forbidden} If $G$ is a contact $B_0$-VPG graph, then
    $G$ is $\{K_5,K_{3,3}, H_0, K_4^-\}$-free.
\end{lem}

\begin{pf}Let $G$ be a contact $B_0$-VPG graph. It immediately follows from Lemma~\ref{o-clique} that $G$ is $K_5$-free.

    Now consider the graph $K_{3,3}$. Let $C$ be a cycle of length
    four in $K_{3,3}$ induced by the vertices $a,b,c,d$. If $K_{3,3}$
    is contact $B_0$-VPG, then, by Lemma~\ref{lem:c4}, in any contact $B_0$-VPG representation
    of $C$,  the union of the paths representing vertices in $C$ must
    enclose a rectangle of the grid. Assume that $P_a,P_c$ are
    horizontal paths, and $P_b,P_d$ are vertical paths. Now, consider
    vertices $e$ and $f$ in $K_{3,3}$ with $e$ being adjacent to $a$
    and $c$, and $f$ being adjacent to $b$ and $d$. Each of the paths
    $P_e,P_f$ must intersect opposite paths of the rectangle. Clearly,
    $P_e$ must be a vertical path and $P_f$ must be a horizontal. If
    $P_e$ is contained inside the rectangle, then it is impossible for
    $P_f$ to intersect $P_b,P_d$ while being inside the rectangle
    without crossing $P_e$. So $P_f$ must be outside the rectangle,
    but then it cannot intersect $P_e$. If $P_e$ lies outside the
    rectangle, then of course $P_f$ has to lie outside the rectangle
    as well, otherwise it cannot intersect $P_e$. But now it cannot
    intersect both $P_b,P_d$ without crossing at least one of them. So
    we conclude that $K_{3,3}$ is not $B_0$-VPG.

    Now let $v,w$ be two adjacent vertices in $G$. Then, in any contact $B_0$-VPG representation of $G$, $P_v$ and $P_w$
    intersect at a grid-point $P$. Clearly, every common neighbour of $v$ and $w$ must also contain $P$. Hence, $v$ and $w$
    cannot have two common neighbours that are non-adjacent. So, $G$ is $K_4^-$-free.

    Finally, consider the graph $H_0$ which consists of three cliques
    of size four, say $A$, $B$ and $C$, with a common vertex $x$.
    Suppose that $H_0$ is contact $B_0$-VPG. Then, it follows from
    Remark~\ref{o-k4} that every vertex in $H_0$ is an end in any
    contact $B_0$-VPG representation of $H_0$. In particular, vertex
    $x$ is an end in any contact $B_0$-VPG representation of $A$, $B$
    and $C$. In other words, the grid-point representing the
    intersection of the paths of each of these three cliques
    corresponds to an endpoint of $P_x$. Since these cliques
    have only vertex $x$ in common, these grid-points are all distinct. But this is a contradiction, since $P_x$ has only
    two endpoints. So we conclude that $H_0$ is not contact $B_0$-VPG, and hence the result follows. \qed   \end{pf}


\section{Chordal graphs}
\label{sec:chordal}


In this section, we will consider chordal graphs and characterise
those that are contact $B_0$-VPG. First, let us point out the
following corollary.

\begin{cor}
    \label{o-block}
    A chordal contact $B_0$-VPG graph is a block graph.
\end{cor}

This follows directly from Lemma~\ref{lem:forbidden} and the definition of block graphs.\\

The following lemma states an important property of minimal chordal non contact $B_0$-VPG graphs that contain neither $K_5$ nor $K_4^-$.

\begin{lem}
    \label{lem:minimal} Let $G$ be a $\{K_5,K_4^-\}$-free
    graph. If $G$ is a minimal non contact $B_0$-VPG graph, then every
    simplicial vertex of $G$ has degree exactly three.
\end{lem}

\begin{pf} Since $G$ is $K_5$-free, every clique in $G$ has size at most
    four. Therefore, every simplicial vertex has degree at most three.
    Let $v$ be a simplicial vertex of $G$. Assume first that $v$ has
    degree one and consider a contact $B_0$-VPG representation of
    $G-v$ (which exists since $G$ is minimal non contact $B_0$-VPG).
    Let $w$ be the unique neighbour of $v$ in $G$. Without loss of
    generality, we may assume that the path $P_w$ lies on some row of
    the grid. Now clearly, we can add one extra column to the grid
    between any two consecutive vertices of the grid belonging to
    $P_w$ and adapt all paths without changing the intersections
    (if the new column is added between column $y_i$ and $y_{i+1}$,
    we extend all paths containing a grid-edge with endpoints in column $y_i$
    and $y_{i+1}$ in such a way that they contain the new edges in the same row
    and between column $y_i$ and $y_{i+2}$
    of the new grid, and any other path remains the same).
    But then we may add a path representing $v$ on this column which only
    intersects $P_w$ (adding a row to the grid and adapting the paths
    again, if necessary) and thus, we obtain a contact $B_0$-VPG
    representation of $G$, a contradiction. So suppose now that $v$
    has degree two, and again consider a contact $B_0$-VPG
    representation of $G-v$. Let $w_1,w_2$ be the two neighbours of
    $v$ in $G$. Then, $w_1,w_2$ do not have any other common neighbour
    since $G$ is $K_4^-$-free. Let $P$ be the grid-point
    corresponding to the intersection of the paths $P_{w_1}$ and
    $P_{w_2}$. Since these paths do not cross and since $w_1,w_2$ do
    not have any other common neighbour (except $v$), there is at
    least one grid-edge having $P$ as one of its endpoints and which
    is not used by any path of the representation. But then we may add
    a path representing $v$ by using only this particular grid-edge
    (or adding a row/column to the grid that subdivides this edge and
    adapting the paths, if the other endpoint of the grid-edge belongs
    to a path in the representation). Thus, we obtain a contact
    $B_0$-VPG representation of $G$, a contradiction. We conclude
    therefore that $v$ has degree exactly three. \qed
\end{pf}

Let $v$ be a vertex of a contact $B_0$-VPG graph $G$. An endpoint
of its corresponding path $P_v$ is \emph{free} in a representation
of $G$, if $P_v$ does not intersect any other path at that
endpoint; $v$ is called \emph{internal} if no representation of
$G$ with a free endpoint of $P_v$ exists. If in a representation
of $G$ a path $P_v$ intersects a path $P_w$ but not at an endpoint
of $P_w$, $v$ is called a \emph{middle neighbour} of $w$.

In the following two lemmas we associate the fact of being or not
an internal vertex of $G$ with the contact $B_0$-VPG
representation of $G$.

\begin{lem}
    \label{lem:endfree} Let $G$ be a chordal contact $B_0$-VPG graph
    and let $v$ be a non internal vertex in $G$. Then, there exists a
    contact $B_0$-VPG representation of $G$ in which all the paths
    representing
    vertices in $G-v$ lie to the left of a free endpoint of $P_v$ (by considering $P_v$ as a horizontal path).
\end{lem}

\begin{figure}
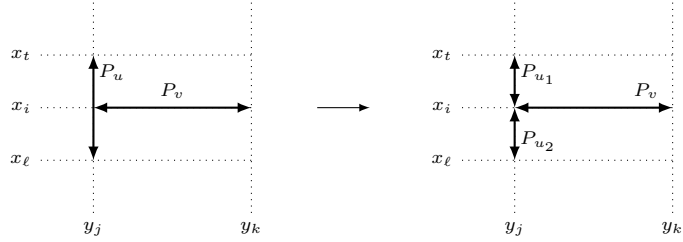

    \begin{center}
    \drawingSplit
        \caption{How to split $P_u$ into two paths.} \label{f-internal}
    \end{center}
\end{figure}

\begin{pf} We will do a proof by induction on the number of vertices of $G$.
    If there is only one vertex in $G$ the result is trivial. Suppose
    $G$ is a graph with at least two vertices. Consider a contact
    $B_0$-VPG representation of $G$. Without loss of generality, we
    may assume that $P_v$ lies on a row $x_i$ between columns
    $y_j,y_k$, $j<k$, and its right endpoint is free. Such a representation exists, since $v$ is not internal.

    If $v$
    is a middle neighbour of another vertex, say $u$, we do the
    following. Assume $P_u$ lies on column $y_j$ between rows
    $x_{\ell}$ and $x_t$, $\ell< t$. We split $P_u$ into two paths,
    $P_{u_1}, P_{u_2}$, such that $P_{u_1}$ goes from row $x_i$ to row
    $x_t$ and $P_{u_2}$ goes from row $x_{\ell}$ to row $x_i$ (see
    Figure~\ref{f-internal}). We denote the corresponding graph by
    $G^*$. If $v$ is not a middle vertex of another vertex,
    then we simply set $G^* = G$.\\

   \noindent  \textbf{Claim.} \emph{The graph $G^*$ is chordal.}\\

    In the second case, it is trivial. In the first case, suppose
    $G^*$ contains a chordless cycle $C$ of length at least $4$. Since
    $G$ is chordal, $C$ contains at least one of $u_1, u_2$.
    Suppose first it contains both $u_1$ and $u_2$. As they are
    adjacent in $G^*$, and contracting them into the vertex $u$
    yields an induced subgraph of $G$, it follows that $C$ has
    length $4$. As in the proof of Lemma~\ref{lem:c4}, it can be
    seen that the paths corresponding to two consecutive vertices in a $C_4$
    cannot be both vertical. So, suppose that $C$ contains only one of $u_1, u_2$, say
    $u_1$. Since $G$ is chordal, $u_2$ has to be adjacent to every
    vertex of $C \setminus N_{G^*}[u_1]$. Since $u_1$ and $u_2$ cannot have two non-adjacent common neighbours, at least one of the neighbours of $u_1$ in $C$ is not adjacent to $u_2$.
Thus, its corresponding path either lies on column $y_j$ having
its lower endpoint in row $x_t$ or lies on some row between
$x_{i+1}$ and $x_t$. In either case, this vertex cannot have a
common neighbour with $u_2$, a contradiction. $\diamondsuit$

    Now, for every vertex $w$ in $N_{G^*}(v)$, consider the connected
    component $C_w$ of $G^* - (N_{G^*}[v] - w)$ containing $w$.
    Notice that $C_w$ is also chordal contact $B_0$-VPG and $w$ is non
    internal in $G^* - (N_{G^*}[v] - w)$. Furthermore, if there are
    two distinct vertices $w$ and $w'$ in $N_{G^*}(v)$, then $C_w$ and $C_{w'}$ are disjoint.
    By contradiction, suppose that a vertex $x$ is in the intersection of $C_{w}$ and $C_{w'}$. Then, there is a path $\alpha_1$ between $w$ and $x$, and a path $\alpha_2$ between $x$ and $w'$. First, suppose $w$ and $w'$ are non adjacent. Joining both paths we can extract a new induced path $\alpha_3$ between $w$ and $w'$ which necessarily has length $\geq 3$. But then, adding $v$ to $\alpha_3$ forms an induced cycle with length $\geq 4$,
     a contradiction. On the other hand, if $w$ and $w'$ are adjacent, first remove the edge $w$ and $w'$. Joining the paths $\alpha_1$ and $\alpha_2$ we can extract an induced path $\alpha_3$ between $w$ and $w'$, which necessarily has length $\geq 4$, since $G$ is $K_4^-$-free
     (see Lemma~\ref{lem:forbidden}) and, therefore, any vertex adjacent to both $w$ and $w'$ must be also adjacent to $v$, implying that it does not belong to $C_w$. Adding the edge between $w$ and $w'$ again, we obtain an induced cycle with length $\geq 4$, a contradiction.

    Finally, considering the case in which $G^* = G$, it is clear that $C_w$ has at least one vertex less than $G$, namely $v$; otherwise, if $u$ was split, the size of $G^*$ is one more than the size of $G$, but then at least two vertices are removed in $C_w$, namely $v$ and one between $u_1$ and $u_2$ (since there is only one vertex in $N_{G^*}[v]$ that we are not removing).

    Then, by induction, there exists a
    contact $B_0$-VPG representation of $C_w$, for each such $w$, with all the paths lying
    to the left of one free endpoint of $P_w$. Now, we replace the
    initial representation of $C_w$ by the new one (the one where all
    the paths lie to the left of one free endpoint of $P_w$) by
    rotating it such that $P_w$ has its free endpoint on the
    grid-point corresponding to the intersection of $P_w$ and $P_v$, and belongs to
    the same side as in the old representation. Notice that we may need to extend the path $P_v$
    to the right before doing the replacement of these new representations to assure
    that they do not overlap.
    Therefore, by extending if necessary the path $P_v$ a little more to the right,
    we obtain a contact $B_0$-VPG representation of $G^*$ in which all the paths lie
    to the left of one free endpoint of $P_v$. In case $P_u$ was split
    into $P_{u_1}$ and $P_{u_2}$, we now glue these two paths together again. \qed\end{pf}

\begin{figure}
    \begin{center}
    \drawingLemmaEndFree
        \caption{Figure illustrating Lemma \ref{lem:endfree}.} \label{f-endfree}
    \end{center}
\end{figure}


\begin{lem}
    \label{lem:internal} Let $G$ be a chordal contact $B_0$-VPG graph.
    A vertex $v$ in $G$ is internal if and only if in every contact
    $B_0$-VPG representation of $G$, each endpoint of the path $P_v$
    either corresponds to the intersection of a representation of
    $K_4$ or intersects a path $P_w$, which represents an internal vertex $w$,
    but not at an endpoint of $P_w$.
\end{lem}

\begin{figure}
    \begin{center}
        \includegraphics[keepaspectratio=true, width=6cm]{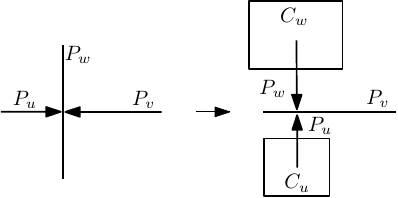}
        \caption{Figure illustrating Lemma \ref{lem:internal}.} \label{f-free1}
    \end{center}
\end{figure}

\begin{pf} The if part is trivial. Assume now that $v$ is an internal vertex
    of $G$ and consider an arbitrary contact $B_0$-VPG representation
    of $G$. Let $P$ be an endpoint of the path $P_v$ and $K$ the maximal clique
    corresponding to all the paths containing the point $P$. Notice that clearly
    $v$ is an end in $K$ by definition of $K$. First, suppose
    there is a vertex $w$ in $K$ which is not an end. Then, it
    follows from Remark~\ref{o-k4} that the size of $K$ is at
    most three. Without loss of generality, we
    may assume that $P_v$ lies on some row and $P_w$ on some column. If $w$ is an internal vertex, we are done. So we may
    assume now that $w$ is not an internal vertex in $G$. Consider $G
    \setminus (K \setminus \{w\})$, and let $C_w$ be the connected
    component of $G \setminus (K \setminus \{w\})$ containing $w$.
    Notice that $w$ is not an internal vertex in $C_w$ either. By
    Lemma~\ref{lem:endfree}, there exists a contact $B_0$-VPG
    representation of $C_w$ with all the paths lying to the left of a free endpoint of
    $P_w$. Now, replace the old representation of $C_w$
    by the new one such that $P$ corresponds to the free endpoint of $P_w$
    in the representation of $C_w$ (it might be necessary to refine --by adding rows and/or columns--
    the grid to ensure that there are no unwanted intersections) and $P_w$ uses the same column as before.
    Finally, if $K$ had size three, say it contains some vertex $u$ in
    addition to $v$ and $w$, then we proceed as follows. Similar to
    the above, there exists a contact $B_0$-VPG representation of
    $C_u$, the connected component of  $G\setminus (K\setminus\{u\})$
    containing $u$, with all the paths lying to the left of a free endpoint of $P_u$,
    since $u$ is clearly not internal in $C_u$. We then replace the old representation of $C_u$ by the new one such that the endpoint of $P_u$ that intersected
    $P_w$ previously corresponds to the grid-point $P$ and $P_u$ lies
    on the same column as $P_w$ (again, we may have to refine the grid). This clearly gives us a contact
    $B_0$-VPG representation of $G$. But now we may extend $P_v$ such
    that it strictly contains the grid-point $P$ and thus, $P_v$ has a
    free endpoint, a contradiction (see Figure~\ref{f-free1}). So $w$ must be an internal vertex.

    Now, assume that all vertices in $K$ are ends.
    If $|K|=4$, we are done. So we may
    assume that $|K|\leq 3$. Hence, there is at least one grid-edge
    containing $P$, which is not used by any paths of the
    representation. Without loss of generality, we may assume that
    this grid-edge belongs to some row $x_i$.  If $P_v$ is horizontal, we may extend it such that it strictly contains $P$. But then $v$ is not internal anymore, a contradiction. If $P_v$ is vertical, then we may extend $P_w$, where $w\in K$ is such that $P_w$ is a horizontal path. But now we are again in the first case discussed above.\qed
\end{pf}

In other words, Lemma~\ref{lem:internal} tells us that a vertex
$v$ is an internal vertex in a chordal contact $B_0$-VPG graph if
and only if we are in one of the following situations:

\begin{itemize}
    \item $v$ is the intersection of two cliques of size four (we say
    that $v$ is of type~1);
    \item $v$ belongs to exactly one clique of
    size four and in every contact $B_0$-VPG representation, $v$ is a
    middle neighbour of some internal vertex (we say that $v$ is of
    type~2);
    \item $v$ does not belong to any clique of size four and
    in every contact $B_0$-VPG representation, $v$ is a middle
    neighbour of two internal vertices (we say that $v$ is of type~3).
\end{itemize}

Notice that two internal vertices of type~1 cannot be adjacent
(except when they belong to a same $K_4$). Furthermore, an
internal vertex of type~1 cannot be the middle-neighbour of some
other vertex.



\begin{figure}[h]
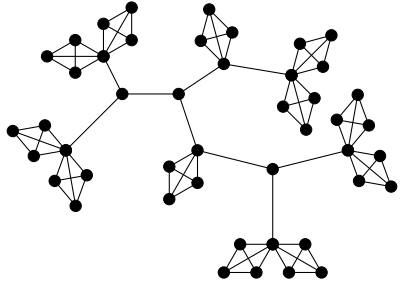

    \begin{center}
\drawingTreeT
\caption{An example of a graph in $\mathcal{T}$.}
        \label{fig:tree}
    \end{center}
\end{figure}

Let $\mathcal{T}$ be the family of graphs containing $H_0$ (see
Figure~\ref{fig:H0}) as well as all graphs that can be partitioned
into a nontrivial tree $T_0$ of maximum degree at most three and
the disjoint union of triangles, in such a way that each triangle
is complete to a vertex $v$ of $T_0$ and anticomplete to $T_0 -
\{v\}$, every leaf $v$ of $T_0$ is complete to exactly two
triangles, every vertex $v$ of degree two in $T_0$ is complete to
exactly one triangle, and vertices of degree three in $T_0$ have
no neighbours outside $T_0$ (see Figure~\ref{fig:tree}).

Notice that all graphs in $\mathcal{T}$ are chordal. We denote by
$B(T)$ the base tree of $T$ in $\mathcal{T}$.

\begin{lem}
    \label{lem:TGraphs} The graphs in $\mathcal{T}$ are not contact
    $B_0$-VPG.
\end{lem}

\begin{pf} By Lemma~\ref{lem:forbidden}, the graph $H_0$ is not contact
    $B_0$-VPG. Consider now a graph $T\in \mathcal{T}$, $T\neq H_0$. Suppose that
    $T$ is contact $B_0$-VPG. Consider an arbitrary contact $B_0$-VPG representation of $T$.
    Consider the base tree $B(T)$ and direct an edge $uv$ of it from $u$ to $v$ if the path $P_v$ contains an endpoint of the path $P_u$ (this way some edges might be directed both ways). If a vertex $v$ has degree $d_B(v)$ in $B(T)$, then by definition of the family $\mathcal{T}$, $v$ belongs to $3-d_B(v)$ $K_4$'s in $T$. Notice that $P_v$ spends one endpoint in each of these $K_4$'s. Thus, any vertex $v$ in $B(T)$ has at most $2-(3-d_B(v))=d_B(v)-1$ outgoing edges. This implies that the sum of out-degrees in $B(T)$ is at most $\sum_{v\in B(T)}(d_B(v)-1)=n-2$, where $n$ is the number of vertices in $B(T)$. But this is clearly impossible since there are $n-1$ edges in $B(T)$ and all edges are directed.\qed
\end{pf}

We will show now how to construct new graphs in $\mathcal{T}$ from others.

\begin{lem}
    \label{lem:TCons}
    \begin{itemize}

      \item[i)] Given $T \in \mathcal{T}$ and $v \in B(T)$ such that $v$ belongs to at least one $K_4$, say $K$,
then the graph $T'$ constructed by removing the other vertices in
$K$ (different from $v$) and adding one vertex $w$ to $B(T)$,
belonging to two copies of $K_4$ (sharing vertex $w$), and
adjacent to $v$, belongs to $\mathcal{T}$.

      \item[ii)] Given $T_1, T_2 \in \mathcal{T}$, $v_1 \in B(T_1)$ and $v_2 \in B(T_2)$ such that $v_1$ and $v_2$ belong to at least one $K_4$ each, say $K_1$ and $K_2$,
then the graph $T'$ constructed by removing the other vertices in
$K_1$ and $K_2$ (different from $v_1$ and $v_2$) and adding one
vertex $w$ to $B(T_1) \cup B(T_2)$, belonging to a $K_4$, and
adjacent to both $v_1$ and $v_2$, belongs to $\mathcal{T}$.
    \end{itemize}
\end{lem}

\begin{pf} \begin{itemize}

  \item[i)] In this case we have $B(T') = B(T) \cup \{w\}$. It is clear that every vertex in $B(T')$ has degree $3$ or less, since we only changed the degree of $v$, which is one less, and the degree of $w$ is one (only adjacent to $v$ in $B(T')$). Moreover, $w$ is a leaf in $B(T')$ and, by construction, it belongs to two
  copies of $K_4$ (sharing vertex $w$). Finally, notice that $v$ has degree $1$ or $2$ in $B(T)$ since vertices of degree $3$ in $B(T)$ does not belong to any $K_4$.
  If $v$ is a leaf in $B(T)$, then $v$ is a degree $2$ vertex in $B(T')$ and, since we removed the other vertices in $K$, it belongs to only one $K_4$ in $T'$. Otherwise, $v$ has degree $2$ in $B(T)$ and therefore, it has degree $3$ in $B(T')$ and does not belong to any $K_4$ in $T'$. Thus, $T' \in \mathcal{T}$.
  \item[ii)] In this case we have $B(T') = B(T_1) \cup B(T_2) \cup \{w\}$. The proof follows in the same manner as the previous
  item.\qed
\end{itemize}
\end{pf}

For the next lemma we need to consider an orientation of some
edges related to a contact $B_0$-VPG representation of $G$, given
by the following rule. If $v, w \in G$ and $v$ is a middle
neighbour of $w$, then we give the orientation from $v$ to $w$.
Let $C_v$ be the reachable vertices starting from $v$, including
$v$. Notice that if $v$ is internal, $C_v = \{v\}$ if and only if
$v$ is of type~1. Also notice that $C_v$ is independent of the
representation for internal vertices. As a consequence of the
previous lemma, we can prove the following.

\begin{lem}
    \label{lem:BaseT} Let $G$ be a chordal contact $B_0$-VPG graph.
    If a vertex $v$ in $G$ is internal, the graph $G'$ constructed by adding a $K_4$, say $K$, containing $v$ to $G$
    contains an induced subgraph $T \in \mathcal{T}$. Moreover, $B(T) = C_v$.

\end{lem}

\begin{pf} We will prove this by induction in the number of vertices in $C_v$.
By Lemma~\ref{lem:internal}, $v$ must be of type~1, 2 or 3. As
noted before, the base case is when $v$ is of type~1. But then $v$
is the intersection of three cliques of size 4 in $G'$, namely $K$
and the two cliques in which $v$ is an end; and thus, $G'$
contains $T = H_0$. Therefore $B(T) = \{v\} = C_v$.

Now, if $v$ is of type~2, $v$ is a middle neighbour of exactly one
other internal vertex $w$. Therefore $C_v = C_w \cup \{v\}$.
Define $G_w$ as the induced subgraph of the connected component of
$G-v$ containing $w$. Notice that $w$ is still internal in $G_w$
since $v$ is a middle neighbour of $w$ in $G$. Then, adding a $K'
= K_4$ containing $w$ to $G_w$ we obtain a $T_1 \in \mathcal{T}$
induced in $G_w$ (and, therefore, also induced in $G$) with
$B(T_1) = C_w$, by inductive hypothesis applied to $w$ in $G_w$.
By Lemma~\ref{lem:TCons} i), we can construct $T \in \mathcal{T}$
by removing the other vertices in $K'$ (different from $w$) and
adding the vertex $v$ (in $G$) to $B(T_1)$, which belongs to two
copies of $K_4$ (one is $K$ and the other is the one in which $v$
is an end), and is adjacent to $w$. Then, $T$ is an induced
subgraph of $G$ and we have $B(T) = C_w \cup \{v\} = C_v$.
Finally, if $v$ is of type~3, $v$ is a middle neighbour of exactly
two other internal vertices $w_1$ and $w_2$. The proof continues
in the same manner as before, applying inductive hypothesis to the
corresponding $G_{w_1}$ and $G_{w_2}$ and then using the second
item of Lemma~\ref{lem:TCons}.\qed
\end{pf}

Using Lemmas~\ref{lem:minimal}--\ref{lem:BaseT}, we are able to
prove the following theorem, which provides a minimal forbidden
induced subgraph characterisation of chordal contact $B_0$-VPG graphs.

\begin{thm}
    \label{thm:contact_chordal}
    Let $G$ be a chordal graph. Let $\mathcal{F}$ = $\mathcal{T} \cup \{K_5,K_4^-\}$. Then, $G$ is a contact $B_0$-VPG graph if and only if $G$ is $\mathcal{F}$-free.
\end{thm}

\begin{pf} Suppose that $G$ is a chordal contact $B_0$-VPG graph. It follows
    from Lemma~\ref{lem:forbidden} and Lemma~\ref{lem:TGraphs} that
    $G$ is $\mathcal{T}$-free and contains neither a $K_4^-$ nor a $K_5$.

    Conversely, suppose now that $G$ is chordal and
    $\mathcal{F}$-free. By contradiction, suppose that $G$ is not
    contact $B_0$-VPG and assume furthermore that $G$ is a minimal non
    contact $B_0$-VPG graph. Let $v$ be a simplicial vertex of $G$
    ($v$ exists since $G$ is chordal). By Lemma~\ref{lem:minimal}, it
    follows that $v$ has degree three. Consider a contact $B_0$-VPG
    representation of $G-v$ and let $K=\{v_1,v_2,v_3\}$ be the set of
    neighbours of $v$ in $G$. Since $G$ is $K_4^-$-free, it
    follows that any two neighbours of $v$ cannot have a common
    neighbour which is not in $K$. First suppose that all the vertices
    in $K$ are ends in the representation of $G-v$. Thus, there exists
    a grid-edge not used by any path and which has one endpoint
    corresponding to the intersection of the paths $P_{v_1}, P_{v_2},
    P_{v_3}$. But now we may add the path $P_v$ using exactly this
    grid-edge (we may have to add a row/column to the grid that subdivides this
    grid-edge and adapt the paths, if the other endpoint of the
    grid-edge belongs to a path in the representation). Hence, we obtain a
    contact $B_0$-VPG representation of $G$, a contradiction.

    Thus, we may assume now that there exists a vertex in $K$ which is
    not an end, say $v_1$. Notice that $v_1$ must be an internal
    vertex. If not, there is a contact $B_0$-VPG representation of
    $G-v$ in which $v_1$ has a free end. Then, using similar arguments
    as in the proof of Lemma~\ref{lem:internal}, we may obtain a
    representation of $G-v$ in which all vertices of $K$ are ends. As
    described previously, we can add $P_v$ to obtain a contact
    $B_0$-VPG representation of $G$, a contradiction. Now, consider the graph $G-K$.
    This graph is clearly chordal contact $B_0$-VPG as being an induced subgraph of $G-v$.
    Then, by Lemma~\ref{lem:BaseT},
    adding the clique $K \cup \{v\}$ (containing the internal vertex $v$) to $G-K$ (which gives the graph $G$) contains an induced subgraph $T \in \mathcal{T}$, a contradiction.  \qed
\end{pf}

Interval graphs form a subclass of chordal graphs. They are defined as being chordal graphs not containing any asteroidal triple, i.e. not containing three pairwise non-adjacent vertices such that there exists a path between any two of them avoiding the neighbourhood of the third one. Clearly, any graph in $\mathcal{T}$ for which the base tree has maximum degree three contains an asteroidal triple. On the other hand, $H_0$ and every graph in $\mathcal{T}$ obtained from a base tree of maximum degree at most two are clearly interval graphs. Denote by $\mathcal{T'}$ the family consisting of $H_0$ and the graphs of $\mathcal{T}$ whose base tree has maximum degree at most two. We obtain the following corollary which provides a minimal forbidden induced subgraph characterisation of contact $B_0$-VPG graphs restricted to interval graphs.

\begin{cor}
    \label{cor:contact_interval} Let $G$ be an interval graph and $\mathcal{F'}$ = $\mathcal{T'} \cup \{K_5,K_4^-\}$.Then, $G$
    is a contact $B_0$-VPG graph if and only if $G$ is $\mathcal{F'}$-free.
\end{cor}


\section{Recognition algorithm}
\label{s:algo}

In this section, we will provide a polynomial-time recognition algorithm for chordal contact $B_0$-VPG graphs which is based on the characterisation given in Section \ref{sec:chordal}. This algorithm takes a chordal graph as input and returns YES if the graph is contact $B_0$-VPG and, if not, it returns NO as well as a forbidden induced subgraph.
The main loop (step 7) will try to find a graph $T\in \mathcal{T}$, $T\neq H_0$. For this purpose, some vertices will be marked and some edges will be directed and coloured.
At the beginning all vertices are unmarked and all edges are undirected and uncoloured.
We will first give the pseudo-code of our algorithm and then explain the different steps.\\

\textbf{Input:} a chordal graph $G=(V,E)$;\\
\textbf{Output:} YES, if $G$ is contact $B_0$-VPG; NO and a forbidden induced subgraph, if $G$ is not contact $B_0$-VPG.

\begin{enumerate}
    \item[$1$.] list all maximal cliques in $G$;
    \vspace*{0.1cm}
    \item[$2$.] if some edge belongs to two maximal cliques, return NO
    and $K_4^-$;
    \vspace*{0.1cm}
    \item[$3$.] if a maximal clique
    contains at least five vertices, return NO and $K_5$;
    \vspace*{0.1cm}
    \item[$4$.] label the vertices such that $l(v)$ =
    number of $K_4$'s that $v$ belongs to;
    \vspace*{0.1cm}
    \item[$5$.]
    if for some vertex $v$, $l(v)\geq 3$, return NO and $H_0$;
    \vspace*{0.1cm}
    \item[$6$.] if $l(v)\leq 1$ $\forall v\in
    V\setminus \{w\}$ and $l(w)\leq 2$, return YES;
    \vspace*{0.1cm}
    \item[$7$.] while there exists an unmarked vertex $v$ with
    $2-l(v)$ outgoing arcs incident to it, do
    \begin{itemize}
        \item[$7.1$] mark $v$ as internal;
        \vspace*{0.1cm}
        \item[$7.2$]
        direct the edges that are currently undirected, uncoloured, not
        belonging to a $K_4$, and incident to $v$ towards $v$;
        \vspace*{0.1cm}
        \item[$7.3$] for any two incoming arcs $wv,w'v$
        such that $ww'\in E$, colour $ww'$;
        \vspace*{0.1cm}
    \end{itemize}
    \item[$8$.] if there exists some vertex $v$ with more than
    $2-l(v)$ outgoing arcs, return NO and find $T\in \mathcal{T}$ by running $BFS$ starting with $v$, following the outgoing arcs,
    and adding for each vertex the corresponding $K_4$'s that it belongs to;
    else return YES.
\end{enumerate}

\begin{figure}[h]
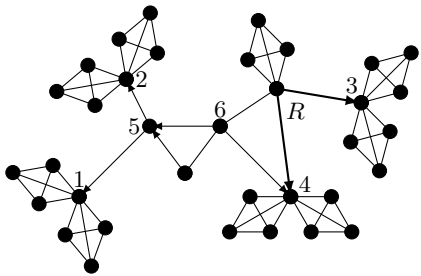

    \begin{center}
    \drawingTreeAlgo
    \caption{
    An example of a possible running of the algorithm. The vertices marked in the algorithm are numbered
    in the order of the marking process. The vertex labeled $R$ corresponds to the root of the tree in the forbidden structure, given in step 8 (whose other vertices are marked as $3$ and $4$).
    } \label{fig-alg}
    \end{center}
\end{figure}

Steps 1-5 can clearly be done in polynomial time (see for example~\cite{G-H-P-chordal} for listing all maximal cliques in a chordal graph). Furthermore, it is obvious to see how to find the forbidden induced subgraph in steps 2, 3 and 5. Notice that if the algorithm has not returned NO after step 5, we know that $G$ is $\{K_4^-,K_5,H_0\}$-free. So we are left with checking whether $G$ contains some graph $T\in \mathcal{T}$, $T\neq H_0$. Since each graph $T\in \mathcal{T}$ contains at least two vertices belonging to two $K_4$'s, it follows that if at most one vertex has label 2, $G$ is $\mathcal{T}$-free (step 6), and thus we conclude by Theorem \ref{thm:contact_chordal} that $G$ is contact $B_0$-VPG.\\

During step 7, we detect those vertices in $G$ that, in case $G$ is contact $B_0$-VPG, must be internal vertices (and mark them as such) and those vertices $w$ that are middle neighbours of internal vertices $v$ (we direct the edges $wv$ from $w$ to $v$). Furthermore, we colour those edges whose endpoints are middle neighbours of a same internal vertex.

Consider a vertex $v$ with $2-l(v)$ outgoing arcs. If a vertex $v$ has $l(v)=2$, then, in case $G$ is contact $B_0$-VPG, $v$ must be an internal vertex (see Lemma \ref{lem:internal}). This implies that any neighbour of $v$, which does not belong to a same $K_4$ as $v$, must be a middle neighbour of $v$. If $l(v)=1$, this means that $v$ belongs to one $K_4$ and is a middle neighbour of some internal vertex. Thus, by Lemma \ref{lem:internal} we know that $v$ is internal. Similarly, if $l(v)=0$, this means that $v$ is a middle neighbour of two distinct internal vertices. Again, by Lemma \ref{lem:internal} we conclude that $v$ is internal. Clearly, step 7 can be run in polynomial time.\\


So we are left with step 8, i.e., we need to show that $G$ is
contact $B_0$-VPG if and only if there exists no vertex with more
than $2-l(v)$ outgoing arcs. First notice that only vertices
marked as internal have incoming arcs. Furthermore, notice that
every maximal clique of size three containing an internal vertex
has two directed edges of the form $wv$, $w'v$ and the third edge
is coloured, where $v$ is the first of the three vertices that was
marked as internal. This is because the graph is $K_4^-$-free and
the edges of a $K_4$ are neither directed nor coloured.

\begin{lem}
    Every vertex marked as internal in step 7 has either label 2 or is the root of a directed induced tree (directed from the root to the leaves) where the root $w$ has degree $2 - l(w)$ and every other vertex $v$ has degree $3 - l(v)$ in that tree, namely one incoming arc and $2 - l(v)$ outgoing arcs.
\end{lem}

\begin{pf}
  By induction in the number of iterations in step 7. In the first iteration, no edge has been directed.
  Therefore, any vertex marked as internal must have label 2, having zero outgoing edges.
  Now assume the result is true for any vertex marked before the $n$-th iteration.
  Let $v$ be the vertex marked in the $n$-th iteration. If $l(v) = 2$ we are done.
  Suppose $l(v) = 1$. Then, there is an outgoing edge from $v$ to a vertex $w$.
  Since only vertices marked as internal have incoming arcs, $w$ must be internal.
  Now, by inductive hypothesis ($w$ was marked in a previous iteration), the result is true for $w$.
  If $l(w) = 2$, $v$ is the root of the tree consisting of the two vertices $v$ and $w$, where $v$ has degree $2 - l(v) = 1$ and $w$ has degree $3 - l(w) = 1$ (one incoming arc). Otherwise, $w$ is the root of a tree $T'$ satisfying the hypothesis of the lemma, but then the tree $T$ constructed from $T'$ by adding $v$ with an outgoing edge to $w$ also clearly satisfies the hypothesis. In a similar manner can be constructed the tree in the case $l(v) = 0$. Finally, let us show that the tree is necessarily induced. Suppose there is an edge not in the tree that joins two vertices of the tree. Since the graph is a block graph, the vertices in the resulting cycle induce a clique, so in particular there is a triangle formed by two edges of the tree and an edge not in the tree. But, as observed above, in every triangle of $G$ having two directed edges, the edges point to the same vertex (and the third edge is coloured, not directed). Since no vertex in the tree has in-degree more than one, this is
  impossible.\qed
\end{pf}

Based on the lemma, it is clear now that if a vertex has more than
$2 - l(v)$ outgoing arcs, then that vertex is the root of a
directed induced tree (directed from the root to the leaves),
where every vertex $v$ has degree $3 - l(v)$, i.e., a tree that is
the base tree $B(T)$ of a graph $T\in \mathcal{T}$. Indeed, notice
that every vertex $v$ in a base tree has degree $3-l(v)$. The fact
that tree is induced can be proved the same way as above. This
base tree can be found by a breadth-first search from a vertex
having out-degree at least $3 - l(v)$, using the directed edges.
Thanks to the labels, representing the number of $K_4$'s a vertex
belongs to, it is then possible to extend the $B(T)$ to an induced
subgraph $T\in \mathcal{T}$. This can clearly be implemented to
run in polynomial time.

To finish the proof that our algorithm is correct, it remains to show that if $G$ contains an induced subgraph in $\mathcal{T}$, then the algorithm will find a vertex with at least $3 - l(v)$ outgoing arcs. This, along with Theorem \ref{thm:contact_chordal}, says that if the algorithm outputs YES then the graph is contact $B_0$-VPG (given that the detection of $K_5$, $K_4^-$ and $H_0$ is clear). Recall that we know that $G$ is a block graph after step 2. Notice that if a block of size 2 in a graph of $\mathcal{T}$ is replaced by a block of size 4, we obtain either $H_0$ or a smaller graph in $\mathcal{T}$ as an induced subgraph. Moreover, adding an edge to a graph of $\mathcal{T}$ in such a way that now contains a triangle, then we obtain a smaller induced graph in $\mathcal{T}$.
Let $G$ be a block graph with no induced $K_5$ or $H_0$. By the remark above, if $G$ contains a graph in $\mathcal{T}$ as induced subgraph, then $G$ contains one, say $T$, such that no edge of the base tree $B(T)$ is contained in a $K_4$ in $G$, and no triangle of $G$ contains two edges of $B(T)$. So, all the edges of $B(T)$ are candidates to be directed or coloured.

In fact, by step 7 of the algorithm, every vertex of $B(T)$ is eventually marked as internal, and every edge incident with it is either directed or coloured, unless the algorithm ends with answer NO before. Notice that by the remark about the maximal cliques of size three and the fact that no triangle of $G$ contains two edges of $B(T)$, if an edge $vw$ of $B(T)$ is coloured, then both $v$ and $w$ have an outgoing arc not belonging to $B(T)$. So, in order to obtain a lower bound on the out-degrees of the vertices of $B(T)$ in $G$, we can consider only the arcs of $B(T)$ and we can consider the coloured edges as bidirected edges. With an argument similar to the one in the proof of Lemma \ref{lem:TGraphs}, at least one vertex has out-degree at least $3 - l(v)$.

\section{Tree-cographs}

In this section, we present a minimal forbidden induced subgraph
characterisation for contact $B_0$-VPG graphs within the class of
tree-cographs.

\emph{Tree-cographs}~\cite{Tinhofer-tree-cographs} are a
generalisation of cographs, i.e. $P_4$-free graphs. They are
defined recursively as follows: trees are tree-cographs; the
disjoint union of tree-cographs is a tree-cograph; and the
complement of a tree-cograph is also a tree-cograph.

It follows from the definition that every tree-cograph is either a
tree, or the complement of a tree, or the disjoint union of
tree-cographs, or the join of tree-cographs. Let us start with the
following two trivial facts.

\begin{fact}
    \label{fact:trees} Every tree is a \conB\/ graph.
\end{fact}

\begin{fact}
    \label{fact:union} The disjoint union of \conB\/ graphs is \conB.
\end{fact}

Now let us consider the complement of trees. We obtain the
following.
\begin{lem}
    \label{lem:co-trees} Let $T$ be a tree. Then $\overline T$ is
    \conB\/ if and only if it is $\{K_5, K_4^-\}$-free.
\end{lem}

\begin{pf}
    If $\overline T$ is
    \conB\/, then it follows from Lemma~\ref{lem:forbidden} that
    $\overline T$ is $\{K_5, K_4^-\}$-free.

    Suppose now that $\overline T$ is $\{K_5,K_4^-\}$-free, then $T$ has
    stability number at most 4. In particular, it has at most four
    leaves. Since it does not have co-($K_4$-e)'s either, we conclude
    that $T$ is either a star with at most $4$ leaves, a $P_4$ or a
    $P_5$. Hence, $\overline T$ is either a $K_4\cup K_1$, a $P_4$ or
    $\overline{P_5}$. Clearly, all these graphs are \conB. \qed \end{pf}

Using the previous results, we are able to obtain the following
characterisation of tree-cographs that are \conB\/.

\begin{thm}
    \label{thm:tree-cographs} Let $G$ be a tree-cograph. Then $G$ is
    \conB\/ if and only if $G$ is $\{K_5, K_{3,3},
    H_0, K_4^-\}$-free.
\end{thm}

\begin{pf} If $G$ is
    \conB\/, then it follows from Lemma~\ref{lem:forbidden} that $G$
    is $\{K_5, K_{3,3}, H_0, K_4^-\}$-free.

    Suppose now that $G$
    is a $\{K_5, K_{3,3}, H_0, K_4^-\}$-free tree cograph on $n$ vertices. We will do a proof by induction on the number of
    vertices of $G$. Let us assume
    the theorem holds for graphs of less than $n$ vertices. If $G$ is
    a tree, the complement of a tree or the disjoint union of
    tree-cographs, then the result holds by Facts~\ref{fact:trees},
    \ref{fact:union}, Lemma~\ref{lem:co-trees} and the induction
    hypothesis. So we may assume now that $G$ is the join of two
    tree-cographs, say $G_1,G_2$.

    Since $G$ is $K_4^-$-free, both $G_1$ and $G_2$ are $P_3$-free,
    i.e., they are the disjoint union of cliques. Furthermore, since
    $G$ is $K_5$-free, it follows that $\omega(G_1)+\omega(G_2) \leq
    4$ and, in particular, none of $G_1,G_2$ contains a $K_4$.

    First suppose that one of $G_1,G_2$, say $G_1$, contains a
    triangle. Then $G_2$ contains no $K_2$. But since $G$ is
    $K_4^-$-free, $G_2$ contains no $2K_1$ either. So $G_2$ is the
    trivial graph. Now, since $G$ is $H_0$-free, $G_1$ contains at
    most two triangles. But then $G$ is clearly  \conB\/. We show in
    Figure~\ref{fig2} how to represent the join of the trivial graph
    and a graph consisting in the disjoint union of at most two
    triangles, an arbitrary number of edges and isolated vertices as a
    \conB\/ graph.

    \begin{figure}[h]
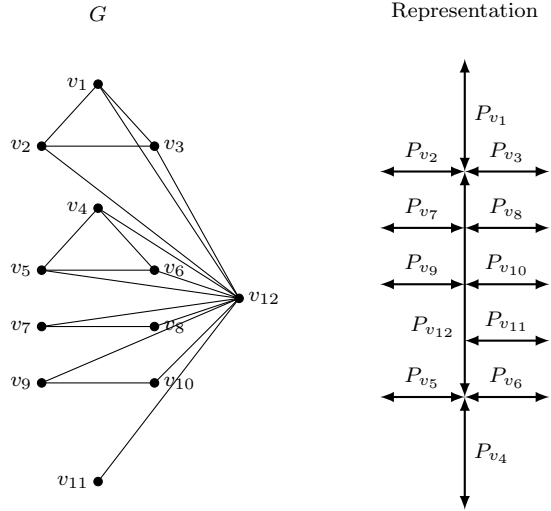

    \begin{center}
\drawingTreeCograph

    \caption{A graph $G$ with $G_1$ with a most two triangles and $G_2=K_1$, and
            a contact $B_0$-VPG representation of $G$.} \label{fig2}
    \end{center}
    \end{figure}

    Next suppose that $\omega(G_1)=\omega(G_2)=2$. Since $G$ is
    $K_4^-$-free, neither $G_1$ nor $G_2$ contains $2K_1$. So $G =
    K_4$, and hence it is \conB.

    Suppose now $\omega(G_1)=2$ and $\omega(G_2)=1$. Since $G$ is
    $K_4^-$-free, $G_2$ contains no $2K_1$, so $G_2$ is the trivial
    graph and hence clearly contact $B_0$-VPG.

%

    Finally, consider the case when  $\omega(G_1)=\omega(G_2)=1$.
    Since $G$ is $K_{3,3}$-free, it follows that $G$ is either the
    star $K_{1,n-1}$ or the complete bipartite graph $K_{2,n-2}$. Thus
    again, $G$ is clearly \conB\/. \qed \end{pf}


From the proofs of the previous results, the following fact can be
deduced.

\begin{cor}
    \label{cor:tree-cographs} Every \conB\/ tree-cograph is the
    disjoint union of trees, $\overline{P_5}$'s, and \conB\/ cographs.
\end{cor}

\section{$P_4$-tidy graphs}


Let $G$ be a graph and let $A$ be a vertex set that induces a
$P_4$ in $G$. A vertex $v$ of $G$ is said to be a \emph{partner}
of $A$ if $G[A\cup\{v\}]$ contains at least two induced $P_4$'s.
The graph $G$ is called \emph{$P_4$-tidy}, if each vertex set $A$
inducing a $P_4$ in $G$ has at most one
partner~\cite{G-R-T-P4tidy}. The class of $P_4$-tidy graphs is an
extension of the class of cographs, i.e. $P_4$-free graphs, and it
contains many other graph classes defined by bounding the number
of $P_4$'s according to different criteria; e.g., $P_4$-sparse
graphs~\cite{Hoang-tesis}, $P_4$-lite graphs~\cite{J-O-p4lite},
and $P_4$-extendible graphs~\cite{J-O-p4extendible}.

A \emph{spider}~\cite{Hoang-tesis} is a graph whose vertex set can be partitioned
into three sets $S$, $C$, and $R$, where $S=\{s_1,\dots,s_k\}$ ($k \geq 2$) is a stable set;
$C=\{c_1,\dots,c_k\}$ is a clique; $s_i$ is adjacent to $c_j$ if
and only if $i=j$ (a \emph{thin spider}), or $s_i$ is adjacent to
$c_j$ if and only if $i \neq j$ (a \emph{thick spider}); $R$ is
allowed to be empty and if it is not, then all the vertices in $R$
are adjacent to all the vertices in $C$ and non-adjacent to all the
vertices in $S$. The triple $(S,C,R)$ is called the \emph{spider
    partition}. By $\thin_k(H)$ and $\thick_k(H)$ we respectively
denote the thin spider and the thick spider with $\vert
C\vert=\vert S\vert=k$ and $H$ the subgraph induced by $R$. If $R$
is an empty set we denote them by $\thin_k$ and $\thick_k$,
respectively. Clearly, the complement of a thin spider is a thick
spider, and vice versa. A \emph{fat spider} is obtained from a
spider by adding a true or false twin of a vertex $v\in S\cup C$.
The following theorem characterises $P_4$-tidy graphs.

\begin{thm}\emph{\cite{G-R-T-P4tidy}}
    \label{thm:P4tidy-char} Let $G$ be a $P_4$-tidy graph with at
    least two vertices. Then, exactly one of the following conditions
    holds:
    \begin{enumerate}
        \item $G$ is disconnected. \item $\overline G$ is disconnected.
        \item $G$ is isomorphic to $P_5$, $\overline{P_5}$, $C_5$, a
        spider, or a fat spider.
    \end{enumerate}
\end{thm}

This allows us to obtain the following characterisation of contact
$B_0$-VPG $P_4$-tidy graphs.

\begin{thm}\label{thm:P4-tidy}
    Let $G$ be a $P_4$-tidy graph. Then $G$ is contact $B_0$-VPG if and only if
    $G$ is $\{K_5,K_{3,3},H_0, K_4^-\}$-free.
\end{thm}

\begin{pf} If $G$ is a contact $B_0$-VPG graph, then it follows from Lemma~\ref{lem:forbidden} that $G$
    is $\{K_5, K_{3,3}, H_0, K_4^-\}$-free.

    Suppose that  $G$ is a $\{K_5,K_{3,3},H_0, K_4^-\}$-free
    $P_4$-tidy graph on $n$ vertices. We will do a proof by induction on the number of
    vertices of $G$. Let us assume the theorem holds
    for graphs of less than $n$ vertices. It follows from
    Theorem~\ref{thm:P4tidy-char} that $G$ is (i) either disconnected;
    (ii) or $\overline G$ is disconnected; (iii) or $G$ is isomorphic
    to $P_5$, $\overline{P_5}$, $C_5$, a spider, or a fat spider.

    If $G$ is disconnected, $G$ is the union of $P_4$-tidy graphs.
    Thus the result holds by Fact~\ref{fact:union} and the induction
    hypothesis.

    If $\overline G$ is disconnected, it follows that $G$ is the join
    of two $P_4$-tidy graphs, say $G_1,G_2$. Then we do exactly the
    same case analysis as in the proof of
    Theorem~\ref{thm:tree-cographs}.

    Now suppose that $G$ is a spider with partition $(C,S,R)$. Since
    $G$ is $K_4^-$-free, $G$ is necessarily a thin spider.
    Furthermore, since $G$ is $K_5$-free, we have $|C| \leq 4$. If
    $|C| = 4$, then $R$ must be empty. If $|C| = 3$, then $|R| \leq 1$
    because $G$ is $\{K_5, K_4^-\}$-free. If $|C| = 2$, then, for
    the same reasons, $|R|\leq 2$ and if $|R|=2$, then $R$ induces
    $K_2$. Notice that for all these cases, the graph obtained is an
    induced subgraph of the graph corresponding to the case $|C| = 4$
    and $R=\emptyset$. We provide a contact $B_0$-VPG representation
    of that case in Figure~\ref{fig5}.

    \begin{figure}[h]
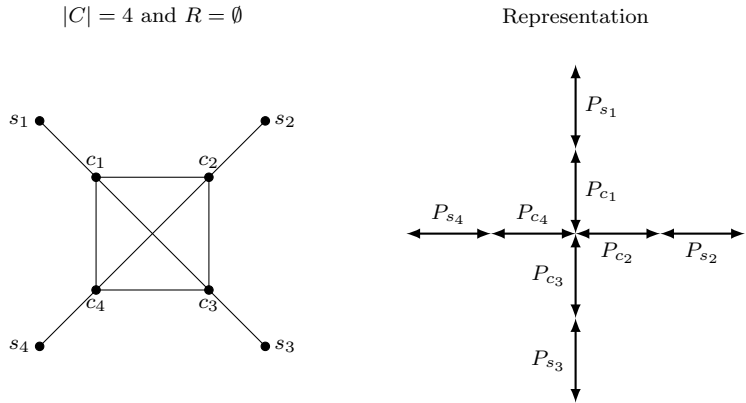

    \begin{center}
\drawingSpider
        \caption{Representation of a thin spider $(C,S,R)$ with $|C| = 4$
            and $R$ empty.} \label{fig5}
    \end{center}
    \end{figure}

    Suppose now that $G$ is a fat spider arising from the thin spider
    with partition $(C,S,R)$. Since $G$ is $K_4^-$-free, it does not
    arise from adding a true twin to a vertex of $C$. For the same
    reason, if $|C| \geq 3$, $G$ does not arise from adding a false
    twin to a vertex of $C$, and if $|C| = 2$, we may add a false twin
    of a vertex of $C$ only if $R$ is empty. We provide a contact
    $B_0$-VPG representation for each of these remaining cases in
    Figure~\ref{fig5a}.

    \begin{figure}[h]
    \begin{center}
{\includegraphics[width=\textwidth]{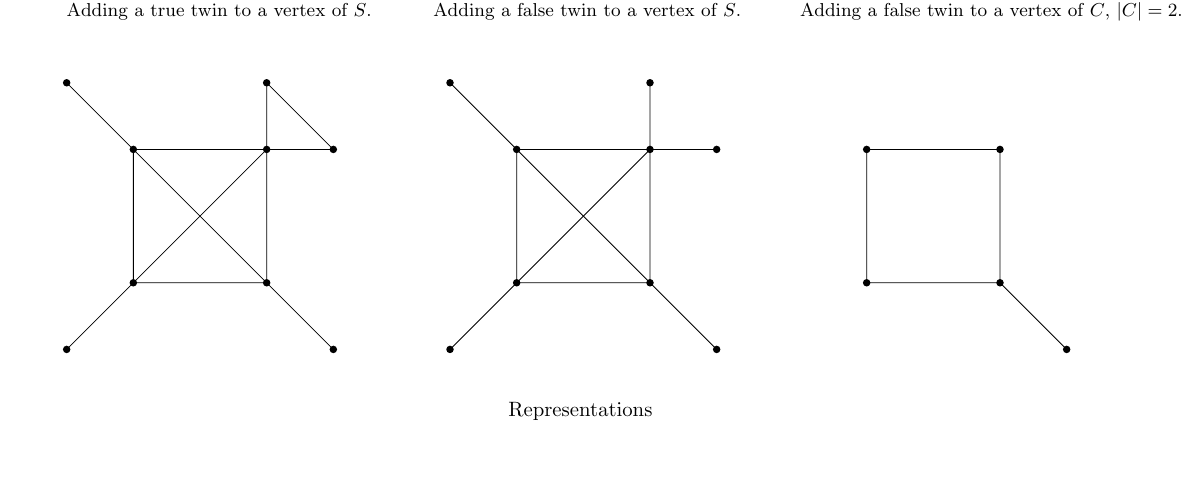}}\\
\drawingFatSpidersRep
        \caption{$G$ is a fat spider arising from the thin spider
            $(C,S,R)$.}\label{fig5a}
\end{center}
    \end{figure}

    Finally, it is easy to see that $P_5$, $\overline{P_5}$, and $C_5$
    are all contact $B_0$-VPG graphs. \qed \end{pf}

    For $P_4$-tidy graphs a linear time recognition algorithm is known~\cite{G-R-T-P4tidy}.
    Using the decomposition properties of the class, the
    characterisation of the possible cases in the proof of
    Theorem~\ref{thm:tree-cographs} for graphs with disconnected
    complement, and the possible cases in the proof of
    Theorem~\ref{thm:P4-tidy} for spiders and fat spiders, we can obtain a linear-time algorithm
    to determine whether a $P_4$-tidy graph is contact $B_0$-VPG.
    Moreover, we can output a minimal forbidden induced subgraph
    in the case the answer is no.


\section{$P_5$-free contact $B_0$-VPG graphs}
\label{s-P5}

In this section, we will present a characterisation of $P_5$-free
contact $B_0$-VPG graphs. Notice that every $P_k$-free graph, with
$1\leq k\leq 2$, is clearly contact $B_0$-VPG. Moreover, a
$P_3$-free graph $G$ is a disjoint union of cliques, therefore $G$
is contact $B_0$-VPG if and only if $G$ is $K_5$-free.

Concerning $P_4$-free graphs, we have the following corollary of
Theorem~\ref{thm:tree-cographs} or Theorem~\ref{thm:P4-tidy},
since $P_4$-free graphs form a subclass of tree-cographs and $P_4$-tidy
graphs.

\begin{thm}
\label{t-P4} Let $G$ be a $P_4$-free graph. Then $G$ is contact
$B_0$-VPG if and only if $G$ is $\{K_5, K_{3,3},
H_0,K_4^-\}$-free.
\end{thm}

Thus, the next graph class to consider is the class of $P_5$-free
graphs. As we will see, the characterisation of $P_5$-free contact
$B_0$-VPG graphs is much more complex than the characterisation of
$P_k$-free graphs, $k\leq 4$. Consider a $P_5$-free graph $G$. If
$G$ is chordal, we obtain a characterisation using
Theorem~\ref{thm:contact_chordal}. Hence, we may assume that $G$
is non chordal. Since $G$ is $P_5$-free it follows that $G$
contains an induced cycle of length $\ell\in \{4,5\}$. In what
follows, we will first analyse the case when $G$ contains an
induced cycle of length four, but no induced cycle of length five.

\begin{lem}
\label{l-main_c4} Let $G$ be a non chordal $\{P_5, C_5, K_{3,3},
K_4^-\}$-free graph. Then, there exists an induced cycle $C$ of
length four in $G$ such that $N[C] = G$.
\end{lem}

\begin{pf}
Since $G$ is not chordal but $\{P_5,C_5\}$-free, it follows that
$G$ must contain an induced cycle of length four. Let $C_0$ be
such a cycle induced by the vertices $v_1,v_2,v_3,v_4$. If $N[C_0]
= G$, we are done. Suppose there exists a vertex $v$ at distance
two of $C_0$. So we may assume, without loss of generality, that
there is a vertex $a$ adjacent to $v_1$ and $v$. It follows from
Remark~\ref{o-Ck} that $a$ must be non-adjacent to at least one of
$v_2,v_4$. Without loss of generality, we may assume that $a$ is
non-adjacent to $v_4$. But then $a$ must be adjacent to $v_3$,
otherwise $v,a,v_1,v_4,v_3$ induce a $P_5$, a contradiction. Thus,
by Remark~\ref{o-Ck}, $a$ is non-adjacent to $v_2$.

Now, consider the cycle $C_1$ induced by the vertices
$a,v_1,v_2,v_3$. If $N[C_1] = G$, we are done. Suppose there is a
vertex $w$ at distance two of $C_1$. Notice that $v,a,v_1,v_4$
induce a $P_4$. Thus, $w$ cannot be adjacent to any of $v,v_4$
otherwise we obtain a $P_5$ or a $C_5$, a contradiction. Hence,
there exists a vertex $b\neq v,v_4$ adjacent to $w$ and to some
vertex in $C_1$. If $b$ is adjacent to exactly one vertex in $C_1$
or to exactly two consecutive vertices in $C_1$, we clearly obtain
a $P_5$, a contradiction. Thus, it follows from Remark~\ref{o-Ck},
that $b$ is adjacent to two nonconsecutive vertices in $C_1$. We
distinguish two cases:

\begin{itemize}
\item[(a)] \textit{$b$ is adjacent to $a$ and $v_2$.} Then $b$
must be adjacent to $v_4$, otherwise $w,b,a,v_1,v_4$ induce a
$P_5$, a contradiction. But now $v_1,v_3,b,a,v_2,v_4$ induce a
$K_{3,3}$, a contradiction.

\item[(b)] \textit{$b$ is adjacent to $v_1$ and $v_3$.} Then $b$
must be adjacent to $v$, otherwise $w,b,v_1,a,v$ induce a $P_5$, a
contradiction. Now consider the cycle $C$ induced by
$a,v_1,b,v_3$. We claim that $N[C] = G$. Suppose there is a vertex
$z$ at distance two of $C$. Then, following the same reasoning as
above, $z$ cannot be adjacent to any of $v_4,v,w,v_2$, since
otherwise we obtain a $P_5$ or $C_5$, a contradiction. Thus, as
before for vertex $b$, there exists a vertex $c$ adjacent to $z$
and to two non-adjacent vertices of $C$. If $c$ is adjacent to
$v_1$ and $v_3$, then $c$ must also be adjacent to $v$, otherwise
$z,c,v_3,a,v$ induce a $P_5$. But now $v_1,v_3,v,a,b,c$ induce a
$K_{3,3}$, a contradiction. Using the same arguments, we can show
that if $c$ is adjacent to $a,b$, then it must be adjacent to
$v_2$, and again we obtain an induced $K_{3,3}$, a contradiction.
Thus $z$ does not exist and hence, $G=N[C]$.\qed\end{itemize}\end{pf}

\begin{figure}
\begin{center}
\includegraphics[scale=0.7]{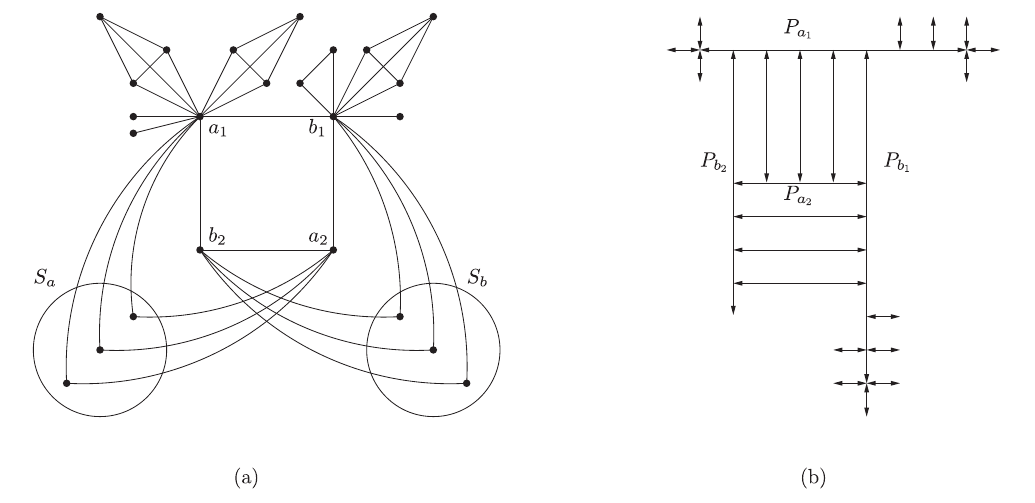}
\caption{(a) An example of a graph from the family
$\mathcal{W}_1$. (b) The corresponding contact $B_0$-VPG
representation.} \label{fig:c4}
\end{center}
\end{figure}

We will define now the following family of graphs.
Start with a cycle $C$ induced by the vertices $a_1,b_1,a_2,b_2$.
Add two (possibly empty) stable sets $S_a$, $S_b$, such that every
vertex in $S_a$ is adjacent to $a_1,a_2$ (but not to $b_1,b_2$),
every vertex in $S_b$ is adjacent to $b_1, b_2$  (but not to
$a_1,a_2$) and $S_a$ is anticomplete to $S_b$. Furthermore, add
two (possibly empty) sets $K_a$, $K_b$ such that $K_a$ (resp.
$K_b$) is complete to $\{a_1\}$ (resp. $\{b_1\}$) and
anticomplete to $\{a_2,b_1,b_2\}$ (resp. $\{a_1,a_2,b_2\}$).
Also, every vertex in $K_a$ (resp. $K_b$) is a simplicial vertex
of degree at most three and $K_a$ (resp. $K_b$) is anticomplete to
$S_a\cup S_b\cup  K_b$ (resp. $S_a\cup S_b\cup K_a$). Finally, add
a (possibly empty) set $K_{ab}$ of vertices forming a clique of
size at most two that is complete to $\{a_1,b_1\}$ and
anticomplete to the rest of the graph. Moreover, neither of $a_1$,
$b_1$ can belong to three cliques of size four and only $a_1$ may
belong to two cliques of size four not containing any vertices
from $K_{ab}$. There are no other edges in the graph. Let us
denote by $\mathcal{W}_1$ the family of graphs described here
before (see Figure~\ref{fig:c4} for an example).

\begin{figure}
\begin{center}
\includegraphics[scale=0.6]{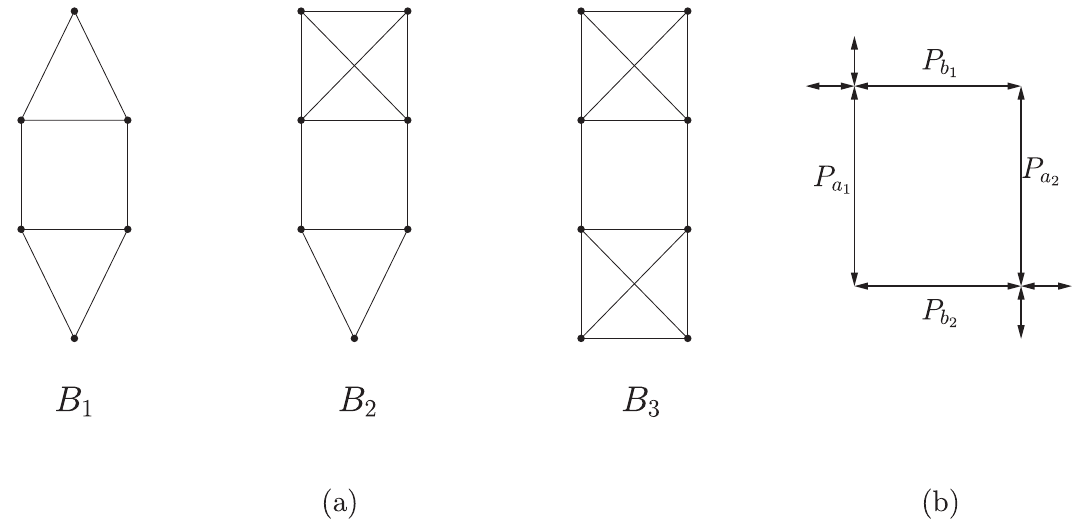}
\caption{(a) The graphs $B_1$, $B_2$ and $B_3$. (b) A contact
$B_0$-VPG representation of $B_3$. } \label{fig:house}
\end{center}
\end{figure}

Let $B_1$, $B_2$ and $B_3$ be the graphs shown in
Figure~\ref{fig:house}. Finally, let $\mathcal{W} = \mathcal{W}_1
\cup \{B_1, B_2, B_3\}$.

\begin{lem}
\label{l-caract_c4} Let $G$ be a non chordal $\{P_5, C_5, K_5,
K_{3,3}, H_0, G_{P_2}, \overline{C_6}, K_4^-\}$-free graph. Then
$G\in \mathcal{W}$.
\end{lem}

\begin{pf}
Let $G$ be a non chordal $\{P_5, C_5, K_5, K_{3,3}, H_0, G_{P_2},
\overline{C_6},K_4^-\}$-free graph. It follows from
Lemma~\ref{l-main_c4} that there exists an induced cycle $C$ of
length four in $G$ such that $N[C] = G$. Let $C$ be induced by
vertices $a_1,b_1,a_2,b_2$.  Let $S_a$ (resp. $S_b$) be the set of
vertices adjacent to $a_1, a_2$ but not $b_1,b_2$  (resp. to
$b_1,b_2$ but not $a_1,a_2$). Notice that $S_a$ (resp. $S_b$) must
be a stable set since $G$ is $K_4^-$-free. Furthermore, $S_a$ is
anticomplete to $S_b$. Indeed, if a vertex $v\in S_a$ is adjacent
to some vertex $w\in S_b$ then $a_1,a_2,w,b_1,b_2,v$ induce a
$K_{3,3}$, a contradiction.

Now, suppose there is a vertex $v$ in $G$ adjacent to only one
vertex in $C$. Without loss of generality, we may assume that $v$
is adjacent to $a_1$. Then, it is not possible to have a vertex
$w\neq v$ in $G$ adjacent only to $a_2$ in $C$, since the vertices
$v,a_1,b_1,a_2,w$ would induce a $P_5$ (in case $v$ and $w$ are
non-adjacent) or a $C_5$ (in case $v$ and $w$ are adjacent).
Therefore, if there is a vertex $w\neq v$ adjacent to only one
vertex in $C$ and different from $a_1$, then we may assume,
without loss of generality, that it is adjacent to $b_1$. Let
$K_a$ (resp. $K_b$) be the set of vertices adjacent to only $a_1$
(resp. $b_1$). If there is a vertex $v \in K_a$ adjacent to a
vertex $w \in K_b$, then $v,w,b_1,a_2,b_2$ induce a $P_5$, a
contradiction. Hence $K_a$ is anticomplete to $K_b$.

Let us now show that all the vertices in $K_a$ are simplicial.
Indeed, suppose that $v \in K_a$ is not simplicial. Then, there
exists $w,u \in N(v)$ such that $u,w$ are non-adjacent. It follows
from the above that $u,w\in K_a$. But then, $v,w,u,a_1$ induce a
$K_4^-$, a contradiction. By symmetry, all vertices in $K_b$ are
simplicial as well. We will distinguish two cases.

First assume now that $G$ is $\overline{P_5}$-free. Thus every
vertex not in $C$ is adjacent to exactly 1 vertex in $C$, since
$G$ is $K_4^-$-free. We claim that $S_a$ is anticomplete to $K_a$.
Indeed, if a vertex $v\in S_a$ is adjacent to some vertex $w\in
K_a$, then $a_1,b_1,a_2,v,w$ induce a $\overline{P_5}$, a
contradiction. Similarly, $S_b$ is anticomplete to $K_b$. Next,
suppose that some vertex $v \in S_a$ is adjacent to some vertex $w
\in K_b$. If $S_b$ is non empty, then for any vertex $u \in S_b$
we obtain a $P_5$ induced by $b_2,u,b_1,w,v$, a contradiction.
Thus, $S_b$ is empty. Then, we may redefine our cycle $C$ by
taking the vertices $a_1,b_1,a_2,v$. Notice that this cycle also
verifies $N[C] = G$. Now, $w \in S_b$ (where $S_b$ is now  the set
of vertices adjacent to $b_1,v$ but not to $a_1,a_2$) and $b_2\in
S_a$. We can proceed similarly if $S_a$ is empty and there are
adjacent vertices in $S_b$ and $K_a$. Now, since $S_b\neq
\emptyset$,  $S_a$ (resp. $S_b$) is anticomplete to $K_b$ (resp.
$K_a$). Since $G$ is $K_5$-free, it follows that the degree of the
simplicial vertices is at most three. Finally, since $G$ is
$\{H_0, G_{P_2}\}$-free, it follows that  only $a_1$ can belong to
two cliques of size four and neither of $a_1$, $b_1$ can belong to
three cliques of size four. Hence, $G\in \mathcal{W}_1$.

Now, suppose that $G$ contains a $\overline{P_5}$ induced by the
cycle $C$ and a vertex $v$ adjacent to $a_1$ and $b_1$. First,
assume there are no other vertices in $G$ adjacent to two
consecutive vertices in $C$. Notice that $v$ cannot be adjacent to
any vertex in $S_a\cup S_b\cup K_a\cup K_b$, since $G$ is
$K_4^-$-free. Moreover, $S_a$ is anticomplete to $K_a$. Indeed, if
$w \in K_a$ is adjacent to $u \in S_a$, then $w,u,a_2,b_1,v$
induce a $P_5$, a contradiction. The same applies to $K_b$ and
$S_b$. Finally, we may assume that $K_a$ (resp. $K_b$) is
anticomplete to $S_b$ (resp. $S_a$) by using the same arguments as
above and redefining the cycle $C$ if necessary. Hence, $G$
belongs to $\mathcal{W}_1$.

Next, assume there is another vertex in $G$ (in addition to $v$)
adjacent to two consecutive vertices in $C$. Notice that $a_1$ and
$b_2$ (resp. $a_2$ and $b_1$) cannot have a common vertex since
$G$ is $P_5$-free. If there is another vertex $w$ adjacent to
$a_1$ and $b_1$, but there is no vertex adjacent to $a_2$ and
$b_2$, then $w$ must be adjacent to $v$, otherwise $a_1,b_1,v,w$
induce a $K_4^-$, a contradiction. Also, $a_1$ (resp. $b_1$)
cannot belong to two cliques of size four whose vertices belong to
$K_a\cup \{a_1\}$ (resp. $K_b\cup \{b_1\}$), since $G$ is
$H_0$-free. Thus, $G$ belongs to $\mathcal{W}_1$, since $G$ is
$K_5$-free and thus no further vertex is adjacent to both $a_1$
and $b_1$. Finally, suppose there is a vertex $w$ adjacent to
$a_2$ and $b_2$. First notice that $v$ and $w$ are non-adjacent,
otherwise $v,w,a_1,b_1,a_2,b_2$ induce a $ \overline{C_6}$, a
contradiction. We claim that all the sets $S_a$, $S_b$, $K_a$ and
$K_b$ must be empty. Indeed, if $u \in S_a$, then $u$ is
non-adjacent to $w$, since $G$ is $K_4^-$-free. But then
$w,a_2,u,a_1,v$ induce a $P_5$, a contradiction. Thus,
$S_a=\emptyset$ and by symmetry we also conclude that
$S_b=\emptyset$. Now suppose $u \in K_a$. Then the vertices $u$,
$a_1$, $b_1$, $a_2$ and $w$ induce a $P_5$ (if $u,w$ are
non-adjacent) or a $C_5$ (if $u,w$ are adjacent). Hence
$K_a=\emptyset$ and by symmetry $K_b=\emptyset$. If there are no
more vertices, $G$ is isomorphic to $B_1$. If there are more
vertices in $G$, then by using the same arguments as before, these
vertices have to be common neighbours of $a_1$ and $b_1$, or $a_2$
and $b_2$. But then $G$ is necessarily isomorphic to either $B_2$
or $B_3$, since $G$ is $\{K_5, \overline{C_6},K_4^-\}$-free (the
same arguments as before apply again).

Finally assume that the $\overline{P_5}$ contained in $G$ is not induced by
the cycle $C$ together with some vertex $v$ adjacent to two
consecutive vertices in $C$. The only possibility is that the
house is induced by $a_1,b_1,a_2,u,w$, with $u\in S_a$ and $w\in
K_a$ (resp. $a_1,b_1,b_2,u,w$, with $u\in S_b$ and $w\in K_b$).
But then, we may redefine our cycle $C$ by taking the vertices
$a_1,b_1,a_2,u$ (resp.$a_1,b_1,b_2,u$). Clearly this new cycle $C$
also verifies that $N[C]=G$. Thus, we can apply the same arguments
as before and show that $G\in \mathcal{W}$.\qed
\end{pf}

\begin{lem}
 \label{l-contact_c4}
Every graph in $\mathcal{W}$ is contact $B_0$-VPG.
\end{lem}

\begin{pf}
Let $G$ be a graph in $\mathcal{W}_1$. We construct a contact
$B_0$-VPG representation of $G$ as follows. First represent the
main cycle $C$ induced by $a_1,b_1,a_2,b_2$: $P_{a_1}$ is a
horizontal path lying on row $x_i$; $P_{a_2}$ is a horizontal path
lying on row $x_j$, $j<i$; $P_{b_1}$ is a vertical path lying on
column $y_k$; $P_{b_2}$ is a vertical path lying on column
$y_{\ell}$, with $\ell>k+|S_a|$; furthermore, we make sure that
$b_1$ and $b_2$ are middle-neighbours of $a_1$ and $a_2$ is a
middle neighbour of $b_1$ and $b_2$; finally the paths $P_{b_1}$
and $P_{b_2}$ use column $y_k$ respectively $y_{\ell}$ down to row
$x_t$ with $t+|S_b|<j$. Now, each vertex in $S_a$ can be
represented by a vertical path on some column $y_r$, with
$k<r<\ell$, and every vertex in $S_b$ can be represented by a
horizontal path on some row $u$ with $t<u<j$. First assume that
$K_{ab}=\emptyset$. Since $P_{a_1}$ has both endpoints free, one
can easily represent two cliques of size four, in case $a_1$
belongs to such cliques and similarly, since $P_{b_1}$ has one
endpoint free, one can easily represent one clique of size four,
in case $b_1$ belongs to such a clique. All other vertices in
$K_a$ or $K_b$ can clearly be represented by extending enough the
paths $P_{a_1}$ and $P_{b_1}$.

Now, assume that $K_{ab}=\{v\}$. Then, given a contact $B_0$-VPG
representation of $G-v$ as described before, we can easily obtain
a contact $B_0$-VPG representation of $G$ as follows: we add a
path $P_v$ lying on column $y_k$ between some row $x_q$ and row
$x_i$, with $i<q$.

Next, assume that $K_{ab}=\{v,v'\}$. Thus, $a_1$ belongs to at
most one clique of size four in $G-\{v,v'\}$ (the vertices of that
clique belong to $K_a$, except for $a_1$).  We obtain a contact
$B_0$-VPG representation as follows. Start with a contact
$B_0$-VPG representation of $G-v'$ as described above. Make sure
that all vertices in $K_a$ are represented by paths intersecting
$P_{a_1}$ to the right of column $y_{\ell}$ (this is clearly
always possible, since $a_1$ belongs to at most one clique of size
four whose vertices (except for $a_1$) belong to $K_a$). Finally,
if necessary, reduce $P_{a_1}$ such that its left endpoint
corresponds to the grid point $(x_i,y_k)$ (this is possible since
$P_{a_1}$ does not intersect any path to the left of that grid
point anymore). Now add $P_w$ as a horizontal path on row $x_i$
with its right endpoint corresponding to the grid point
$(x_i,y_k)$.

Finally, if $G$ is one of the graphs $B_1$, $B_2$ or $B_3$, then
$G$ is clearly contact $B_0$-VPG as can be seen in
Figure~\ref{fig:house}(b). Notice that $B_1,B_2$ are induced
subgraphs of $B_3$.\qed
\end{pf}

From the lemmas above, we conclude the following.

\begin{cor}
\label{coroP4} Let $G$ be a non chordal $\{P_5, C_5, K_5, K_{3,3},
H, G_{P_2}, \overline{C_6}, K_4^-\}$-free graph. Then $G$ is
contact $B_0$-VPG.
\end{cor}

Let us now focus on $P_5$-free graphs containing an induced cycle
of length five.

\begin{lem}
\label{l-c5} Let $G$ be a $\{P_5, K_4^-\}$-free graph. Let $C$ be
an induced cycle of length five in $G$ such that no vertex is
adjacent to exactly three non consecutive vertices in $C$. Then,
$N[C]=G$ and every vertex $v\in N(C)$ is adjacent to exactly two
non-consecutive vertices in $C$.
\end{lem}

\begin{pf}
Let $C$ be induced by $v_1,\cdots,v_5$ and let $v$ be a vertex in
$N(C)$. It follows from Remark~\ref{o-Ck} that $v$ cannot be
adjacent to three consecutive vertices in $C$. If $v$ is adjacent
to exactly one vertex or to two consecutive vertices in $C$, then
we clearly obtain a $P_5$, a contradiction. Thus, $v$ has exactly
two non consecutive neighbours in $C$.

Now assume that there exists a vertex $u$ which is at distance two
of $C$. Thus, there is a vertex $w \in N(C)$ adjacent to $u$ and
to two non consecutive vertices in $C$, say $v_1,v_3$. But then,
$v,w,v_1,u_5,v_4$ induce a $P_5$, a contradiction. Therefore
$N[C]=G$.\qed
\end{pf}

Let $K_{3,3}^*$ be the graph obtained by subdividing exactly one
edge in the graph $K_{3,3}$. We will now define several families of
graphs. Start with a cycle $C$ of length five induced by the
vertices $a,v,b,c,w$. Add two (possibly empty) stable sets $S_v$,
$S_w$ such that $S_v$ is complete to $\{a,b\}$, $S_w$ is complete
to $\{a,c\}$ and $S_v$ is anticomplete to $S_w$. There are no
other edges. Let us denote by $\mathcal{L}_1$ the family of graphs
described here before.

\begin{figure}
\begin{center}
\includegraphics[scale=0.75]{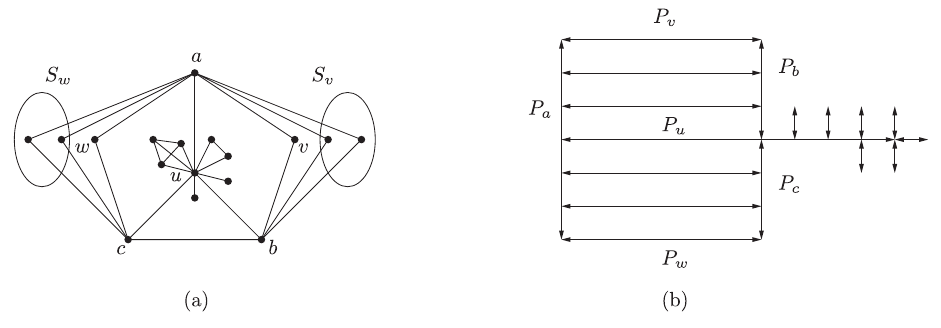}
\caption{(a) An example of a graph in $\mathcal{L}_2$. (b) The
corresponding contact $B_0$-VPG representation. }
\label{fig:c5house}
\end{center}
\end{figure}

Let $G\in \mathcal{L}_1$ and let $G'$ be the graph obtained from
$G$ by adding a vertex $u$ adjacent to $a$, $b$ and $c$.
Furthermore, add a (possible empty) set $K_u$, such that $K_u$ is
complete to $\{u\}$ and anticomplete to $V(C)\cup S_v\cup S_w$.
Also, every vertex in $K_u$ is a simplicial vertex of degree at
most three. Moreover, $u$ can belong to only one clique of size
four. There are no other edges. Let us denote by $\mathcal{L}_2$
the family of graphs described here before (see
Figure~\ref{fig:c5house}(a) for an example).

Next, consider a graph $G'$ in $\mathcal{L}_2$ with $S_v = S_w =
\emptyset$ and $u$ not belonging to any clique of size four. Add a
vertex $z$ adjacent to $v$, $w$ and $u$. There are no other edges.
Let us denote by $\mathcal{L}_3$ the family of graphs obtained
that way and let $\mathcal{L} = \mathcal{L}_1 \cup \mathcal{L}_2
\cup \mathcal{L}_3$.

\begin{figure}
\begin{center}
\includegraphics[scale=0.70]{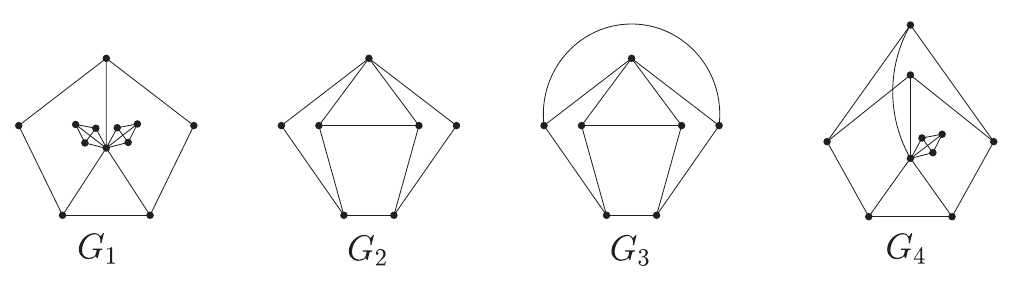}
\caption{The graphs $G_1$, $G_2$, $G_3$ and $G_4$.}
\label{fig:Ggraphs}
\end{center}
\end{figure}

Finally, let $G_1$, $G_2$, $G_3$ and $G_4$ be the graphs shown in
Figure~\ref{fig:Ggraphs}.

\begin{lem}
\label{l-caract_c5} Let $G$ be a $\{P_5, K_5, K_{3,3}^*,
\overline{C_6}, G_1, G_2, G_3, G_4, K_4^-\}$-free graph and assume
$G$ contains a cycle of length five. Then $G \in \mathcal{L}$.
\end{lem}

\begin{pf}
Let $C$ be an induced cycle of length five with vertices
$a,v,b,c,w$. Clearly, no vertex in $N(C)$ is adjacent to exactly
one vertex in $C$ or to two consecutive vertices in $C$, since $G$
is $P_5$-free. Consider first the vertices adjacent to two
non-consecutive vertices in $C$. For any two vertices $u,z$ that
are adjacent to the same two non-consecutive vertices in $C$, we
have that $u$ and $z$ are non-adjacent otherwise we obtain
$K_4^-$, a contradiction. Suppose that there exist vertices $u,z$
such that they have distinct neighbours in $C$, say $u$ is
adjacent to $a$ and $b$, and $z$ is adjacent to $v$ and $c$. If
$u,z$ are adjacent, then together with the vertices of $C$, they
induce a $K_{3,3}^*$, a contradiction. If $u,z$ are non-adjacent,
then $u,a,v,z,c$ induce a $P_5$, a contradiction. Thus, we may
assume now, without loss of generality, that every vertex adjacent
to two nonconsecutive vertices in $C$ is either adjacent to both
$a$ and $b$ or adjacent to both $a$ an $c$. Let $S_v$ (resp.
$S_w$) be the set of vertices not in $C$ adjacent to $a,b$ (and
not to $v,c,w$) (resp. $a,c$ (and not to $b,v,w$)). It follows
from the above that $S_v$ and $S_w$ are stable sets. Finally, if
there is a vertex $u\in S_v$ adjacent to some vertex $z\in S_w$,
then we obtain $G_2$, a contradiction. So $S_v$ is
anticomplete to $S_w$.

First assume that there exists no vertex in $G$ that is adjacent
to three non-consecutive vertices in $C$. It immediately follows
from Lemma~\ref{l-c5} that $G=N[C]$ and that every vertex not in
$C$ is adjacent to two non-consecutive vertices in $C$. Thus, $G
\in \mathcal{L}_1$.

Now, suppose that there exists a vertex $u$ adjacent to three
non-consecutive vertices in $C$, say $a$, $b$ and $c$. We will first
show that there cannot be another vertex adjacent to three
non-consecutive vertices. If there is another vertex $z$ adjacent
to $a$, $b$ and $c$, then $u$ and $z$ must be adjacent otherwise
$u,z,c,b$ induce a $K_4^-$, a contradiction. But now the vertices
$a$, $u$, $z$ and $b$ induce a $K_4^-$, again a contradiction.
Now, suppose $z$ is adjacent to $v$, $b$ and $w$. Then, $z$ and
$u$ are non-adjacent, since otherwise $u,z,b,c$ induce a $K_4^-$, a
contradiction. But now, the vertices of $C$ together with $u$ and
$z$ induce $G_3$, a contradiction as well. By symmetry, we
conclude that $z$ cannot be adjacent to $v$, $w$ and $c$. Finally,
if $z$ is adjacent to $a$, $v$ and $c$, the vertices $a$, $v$,
$z$, $u$, $b$ and $c$ induce a $\overline{C_6}$ if $z$ and $u$ are
non-adjacent, a contradiction. But if $z$ and $u$ are adjacent,
then $u,z,b,c$ induce a $K_4^-$, again a contradiction. By
symmetry, $z$ cannot be adjacent to $a$, $w$ and $b$. Hence, we
conclude that $u$ is the unique vertex adjacent to three
non-consecutive vertices in $C$.

Now we will distinguish several cases, depending on which vertices
$u$ is adjacent to. First, assume that $u$ is adjacent to $a$, $b$
and $c$, and that $S_v\cup S_w$ is non empty. Notice that $u$
cannot be adjacent to any vertex in $S_v\cup S_w$, since $G$ is
$K_4^-$-free. It follows from Remark~\ref{o-Ck} and the fact that
$G$ is $P_5$-free that any vertex in $G$ not belonging to
$V(C)\cup S_v\cup S_w\cup \{u\}$ has to be adjacent to $u$ and
anticomplete to $V(C)\cup S_v\cup S_w$. Let $K_u=N(u)\setminus
V(C)$ be the set of these vertices and consider $z\in K_u$. Then
$z$ is simplicial. Indeed, if $z$ is not simplicial, it follows
that there exist vertices $z',z''\in K_u\cap N(z)$ such that
$z',z''$ are non-adjacent. But then $z,z',z'',u$ induce $K_4^-$, a
contradiction. Furthermore, since $G$ is $K_5$-free, it follows
that every vertex $z\in K_u$ has degree at most three. Finally,
notice that $u$ can only belong to at most one clique of size
four, since $G$ is $G_1$-free. Thus, we conclude that $G \in
\mathcal{L}_2$.

Notice that if $S_v = S_w = \emptyset$, we can relabel the
vertices in $C$ such that $u$ is adjacent to $a$, $b$ and $c$, and
we obtain a graph in $\mathcal{L}_2$ as before. Thus, we may
assume, without loss of generality, that there is a vertex $z \in
S_v$. Now, we will consider different cases:

\begin{itemize}

\item If $u$ is adjacent to $v$, $b$ and $w$, or if $u$ is
adjacent to $a$, $v$ and $c$, then we obtain $G_2$ (notice that
$z$ and $u$ cannot be adjacent since the graph is $K_4^-$-free),a
contradiction.

\item If $u$ is adjacent to $a$, $b$ and $w$, then $S_w =
\emptyset$, otherwise $a,v,b,c,u,w,t$, where $t\in S_w$, induce
$G_2$ a contradiction. Now, we can relabel the vertices in $C$
such that $u$ is adjacent to $a$, $b$ and $c$, without changing
$S_v$, and we obtain a graph in $\mathcal{L}_2$ as before.

\item If $u$ is adjacent to $v$, $c$ and $w$, and $z$ is non-adjacent to $u$, then $z,a,v,u,c$ induce a $P_5$, a contradiction.
So $z$ and $u$ must be adjacent. Notice again that $S_w =
\emptyset$. Indeed, if $t\in S_w$, then $t,a,v,b,u,c,w$ induce
$G_2$, a contradiction. Moreover, $|S_v|=1$: if $z'\in S_v$,
$z\neq z'$, then $z'$ must be adjacent to $u$ as well, but now
$v,z,z',a,b,u$ induce a $K_{3,3}$, a contradiction. So we can
relabel the vertices in $C$ such that $u$ is adjacent to $a$, $b$,
$c$. With this new labeling, $S_v = S_w = \emptyset$ and $z$ is
adjacent to $v$, $w$ and $u$. Clearly, any vertex not belonging to
$V(C)\cup \{u,z\}$ has to be adjacent to $u$, since $G$ is
$P_5$-free. Let $K_u$ be the set of these vertices. Using the same
arguments than above, one can show that ever vertex in $K_u$ is
simplicial and have degree at most three since the graph is
$K_5$-free. Finally, $u$ cannot belong to a clique of size four,
since $G$ is $G_4$-free. So we conclude that $G \in
\mathcal{L}_3$.\qed
\end{itemize}
\end{pf}

\begin{lem}
\label{l-contact_c5} Every graph in $\mathcal{L}$ is contact
$B_0$-VPG.
\end{lem}

\begin{figure}
\begin{center}
\includegraphics[scale=0.70]{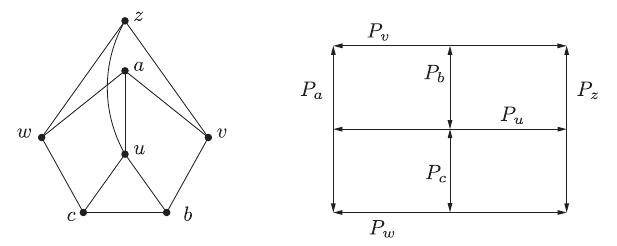}
\caption{A graph in $\mathcal{L}_3$ and the corresponding contact
$B_0$-VPG representation.} \label{fig:L3}
\end{center}
\end{figure}

\begin{pf} Let $G\in \mathcal{L}_1$. We construct a contact $B_0$-VPG
representation of $G$ as follows.  Vertex $b$ is represented by a
path $P_b$ lying on column $y_j$ between rows $x_k$ and $x_t$,
with $t>k+|S_v|$; vertex $c$ is represented by a path $P_c$ lying
on column $y_j$ between rows $x_t$ and $x_{\ell}$, with
$\ell>t+|S_w|$; vertex $a$ is represented by a path $P_a$ lying on
column $y_i$, $i<j$, between rows $x_k$ and $x_{\ell}$; vertex $v$
is represented by a path $P_v$ lying on row $x_k$ between rows
$y_i$ and $y_j$ and vertex $w$ is represented by a path $P_w$
lying on row $x_{\ell}$ between rows $y_i$ and $y_j$. Now each
vertex in $S_v$ is represented by a path between columns $y_i$ and
$y_j$ lying on one of the $|S_v|$ rows between $x_k$ and $x_t$,
and each vertex in $S_w$ is represented by a path between columns
$y_i$ and $y_j$ lying on one of the $|S_w|$ rows between $x_t$ and
$x_{\ell}$.

If $G\in\mathcal{L}_2$, consider a representation of $G-(K_u\cup
\{u\})$ as described above. Now, it is possible to add $P_u$ on
row $x_t$, such that $b$ and $c$ are middle-neighbours of $u$, and
$u$ is a middle-neighbour of $a$. If $u$ belongs to one clique of
size four, then it is possible to represent this clique using the
right endpoint of $P_u$. All the other vertices of $K_u$ can
easily be represented by eventually extending the path $P_u$ to
the right.

Finally, if $G\in \mathcal{L}_3$, consider the contact $B_0$-VPG
representation of the graph shown in Figure~\ref{fig:L3}. Clearly,
it is possible to add the paths representing the vertices of $K_u$,
since $u$ does not belong to any clique of size four.\qed
\end{pf}

\begin{lem}
\label{l-notcontact} The graphs $K_{3,3}^*, \overline{C_6}, G_1,
G_2, G_3, G_4$  are not contact $B_0$-VPG.
\end{lem}

\begin{pf} Consider the graph $K_{3,3}$ with vertices $a,c,e$ on one side of
the bipartition and $b,d,f$ on the other side. Assume that the
edge $ef$ is subdivided to obtain $K_{3,3}^{*}$. Consider the
cycle induced by the vertices $a,b,c,d$. Following the same
approach as in Lemma~\ref{lem:forbidden}, we may assume that
$P_a,P_c$ are horizontal paths, $P_b,P_d$ are vertical paths and
$P_e$ is a horizontal path lying inside the rectangle, and $P_f$
is a vertical path lying outside the rectangle. But now it is
clearly impossible to add a path intersecting $P_e$ and $P_f$
without intersecting any other path. Thus, $K_{3,3}^*$ is not
$B_0$-VPG.

Next consider the graph $\overline{C_6}$ with vertex set
$a,b,c,d,v,w$ such that $a,b,c,d$ induce a cycle of length four,
$v$ is a common vertex of $a$ and $b$, $w$ is a common neighbour
of $c$ and $d$, and $v$ is adjacent to $w$. If $\overline{C_6}$ is
contact $B_0$-VPG, then we may assume that in a contact $B_0$-VPG
representation, the paths $P_a,P_c$ are horizontal and the paths
$P_b,P_d$ are vertical. Since $b,c,v,w$ induce a cycle of length
four, we conclude from the above that $P_v$ has to be horizontal.
But since $a,d,v,w$ induce a cycle of length four as well, we also
conclude that $P_v$ has to be vertical, a contradiction. Hence,
$\overline{C_6}$ is not $B_0$-VPG.

Suppose now that the graph $G_1$ is contact $B_0$-VPG. Without
loss of generality, we may assume that $P_u$ lies on some row
$x_i$. Since $u$ belongs to two cliques of size four, it follows
from Remark~\ref{o-k4} that both endpoints of $P_u$ are not free.
Thus, $a,b$ and $c$ are middle neighbours of $u$, i.e.
$P_a,P_b,P_c$ are necessarily vertical paths. Thus, $P_v,P_w$ must
be horizontal paths, but this is impossible since no two paths can
cross. We conclude that $G_1$ is not contact $B_0$-VPG.

Using similar arguments, we conclude that if $G_4$ is contact
$B_0$-VPG, then $b,c$ have to be middle neighbours of $u$, $u$ has
to be a middle neighbour of $a$ and $P_v,P_w$ have to be
horizontal paths. But now it is clearly impossible to add $P_z$
such that it intersects $P_v,P_w,P_u$ without crossing any path.
Hence, $G_4$ is not contact $B_0$-VPG.

Finally, consider the graphs $G_2,G_3$ and suppose that they are
contact $B_0$-VPG. First consider a contact $B_0$-VPG
representation of $G_2-v$ (resp. $G_3-v$). Since $t$ is adjacent
to three non-consecutive vertices of a induced cycle of length
five, we may assume, without loss of generality, that we have the
following configuration: $P_a,P_c,P_z$ are horizontal paths with
$P_a,P_z$ lying on a same row; $P_b,P_w$ are vertical paths; $P_t$
is a vertical path with one endpoint corresponding to the
endpoints of $P_a,P_z$ that intersect; $t$ is a middle neighbour
of $c$. But now it is clearly impossible to add a path representing
vertex $v$, since it has to intersect $P_a$ and $P_b$. Therefore,
$G_2,G_3$ are not contact $B_0$-VPG. \qed
\end{pf}

We are now ready to prove the main result of this section.

\begin{thm}
\label{t-contact_p5housefree} Let $G$ be a $P_5$-free graph. Let
$\mathcal{G} = \{K_5,H_0,G_{P_2},K_{3,3},K_{3,3}^*,$
$\overline{C_6},G_1,G_2,G_3,G_4, K_4^-\}$. Then $G$ is contact
$B_0$-VPG if and only if $G$ is $\mathcal{G}$-free.
\end{thm}

\begin{pf}
For the only if part, we use Theorem~\ref{thm:contact_chordal},
Lemma~\ref{lem:forbidden} and Lemma~\ref{l-notcontact}.

Suppose now that $G$ is a $P_5$-free graph which is also
$\mathcal{G}$-free. If $G$ is chordal, the result follows from
Theorem~\ref{thm:contact_chordal}, since $G$ is $\mathcal{F}$-free
(indeed, the graphs in $\mathcal{F}$ different from $H_0$ and
$G_{P_2}$ contain an induced $P_5$). Now, assume that $G$ is not
chordal. If $G$ is $C_5$-free, by Corollary~\ref{coroP4}, $G$ is
contact $B_0$-VPG. Similarly, if $G$ contains a $C_5$, by
Lemmas~\ref{l-caract_c5} and~\ref{l-contact_c5}, $G$ is also
contact $B_0$-VPG.\qed
\end{pf}

\section{Conclusions and Future work} \label{s:conclusion}

In this paper, we considered some special graph classes, namely
chordal graphs, tree-cographs, $P_4$-tidy graphs and $P_5$-free
graphs. We gave a characterisation by minimal forbidden induced
subgraphs of those graphs from these families that are contact
$B_0$-VPG. Moreover, we presented a polynomial-time algorithm for
recognising chordal contact $B_0$-VPG graphs based on our
characterisation. Notice that for the other graph classes
considered here, the characterisation immediately yields a
polynomial-time recognition algorithm.

In order to get a better understanding of the structure of general
contact $B_0$-VPG graphs, one way could be to find further
characterisations by forbidden induced subgraphs of contact
$B_0$-VPG graphs within other interesting classes. Since classical
graph problems are difficult in contact $B_0$-VPG graphs (see for
instance~\cite{D-G-M-R-cpg}), these further insights in their
structure may lead to good approximation algorithms for these
problems.

\subsection*{Acknowledgements}
This work was partially supported by ANPCyT PICT-2015-2218 and
UBACyT Grants 20020130100808BA and 20020160100095BA (Argentina).




\end{document}